\documentclass[sn-mathphys,Numbered]{sn-jnl}

\usepackage{manyfoot}%
\usepackage[figuresright]{rotating}
\usepackage{url}
\usepackage{amsmath,amsthm}
\usepackage{amssymb}
\usepackage{multirow,longtable}
\usepackage{subfigure}
\usepackage{fancybox}
\usepackage{booktabs}
\usepackage{color}
\usepackage{xcolor,colortbl}
\usepackage[makeindex]{imakeidx}
\usepackage{natbib}

\newcommand{\gray}{\cellcolor{gray!25}}

\begin{document}

\title{Outsourcing policies for the Facility Location Problem with Bernoulli Demand}

\author*[1]{\fnm{Maria} \sur{Albareda-Sambola}}\email{maria.albareda@upc.edu}

\author*[2]{\fnm{Elena} \sur{Fernández}}\email{elena.fernandez@uca.es}

\author*[3]{\fnm{Francico} \sur{Saldanha-da-Gama}}\email{Francisco.Saldanha-da-Gama@Sheffield.ac.uk}

\affil[1]{\orgdiv{Dept. d'Estad\'{i}stica i Investigaci\'{o} Operativa}, \orgname{Universitat Polit\`{e}cnica de Catalunya}, \orgaddress{\street{Carrer Colom, 11}, \city{Terrassa}, \postcode{08022}, \country{Spain}}}

\affil[2]{\orgdiv{Dept. de Estad\'{i}stica e Investigaci\'{o}n Operativa}, \orgname{Universidad de C\'adiz}, \orgaddress{\street{Pol\'igono R\'io San Pedro}, \city{Puerto Real}, \postcode{11510}, \country{Spain}}}

\affil[3]{\orgdiv{Management School}, \orgname{Sheffield University}, \orgaddress{\street{Conduit Road}, \city{Sheffield}, \postcode{S10 1FL}, \state{England}, \country{United Kingdom}}}

\abstract{This paper focuses on the Facility Location Problem with Bernoulli Demand, a discrete facility location problem with uncertainty  where the joint distribution of the customers' demands is expressed by means of a set of possible scenarios. A two-stage stochastic program with recourse is used to select the facility locations and the \emph{a priori} assignments of customers to open plants, together with the \emph{a posteriori} strategy to apply in those realizations where the \emph{a priori} solution is not feasible.
Four alternative outsourcing policies are studied for the recourse action, and a mathematical programming formulation is presented for each of them. Extensive computational experiments have been carried-out to analyze the performance of each of the formulations and to compare the quality of the solutions produced by each of them relative to the other outsourcing policies.}




\maketitle

\section{Introduction} \label{sect:introduction}
Broadly speaking, facility location problems look for the best locations for a set of facilities that must satisfy service requests of a given set of customers \citep[see, e.g.,][]{LaporteNickelSaldanha:2019}. It is often assumed that customers' demand is part of the input data and, thus, is known in advance.
Nonetheless, in practice, customers' demand is subject to a high level of
uncertainty, so the above assumption very seldom holds. Examples of location problems with non-deterministic demands include any logistic related location problem where demand levels might change over different time periods (postal services, supermarkets, warehouses to distribute goods with seasonal-dependent demand, airports, etc.)  \citet{BrandeauChiu:1989}, \citet{Louveaux:1993}, and \citet{Snyder:2006} have surveyed different aspects of Stochastic Location Problems.

It is easy to find situations, for instance, in logistic applications, where requests of service are unitary in the sense that each service request consumes one resource unit (e.g. one worker) from the service center. This is the case, for example, of delivery services, door-to-door mail services or home assistance services. Facility location problems with unit-demand customers have been studied in the literature motivated by different types of applications namely in Telecommunications \citep[e.g.,][]{Fortz:2015} and Healthcare \citep[e.g.,][]{Ahmadi-JavidSeyediSyam:2017}.
Other work not focusing on a specific application can also be found in \cite{Drexl:2011} and \cite{GouveiaSaldanha-da-Gama:2006}.
When demand is stochastic, then, each possible scenario is characterized by the set of customers with demand, and requests can be modeled by means of binary vectors. In that case, the components of such a vector are Bernoulli random variables. 
Two-stage stochastic models with recourse function \citep[see, e.g.,][]{BirgeLouveaux:2011,HaneveldVanDerVlerkRomeijnders:2020} have been proposed for addressing several stochastic combinatorial optimization problems with Bernoulli demand. This is the case of the probabilistic traveling salesman problem \cite{Jaillet:1988,LaporteLouveauxMercure:1994,BianchiCampbell:2007},  location-routing problems with Bernoulli demands \cite{BermanSimchi-Levi:1988,AlbaredaFernandezLaporte:2007}, or the stochastic generalized assignment Problem \cite{AlbaredaVanDerVlerkFernandez:2006}.

An example of an  stochastic facility location problem with unitary demands is the Facility Location Problem with Bernoulli Demand (FLPBD) that we address in this paper. The FLPBD is a discrete facility location problem where the uncertain service demand of each customer follows a Bernoulli distribution, and the joint distribution of customers' demands is expressed by means of a set of possible scenarios.
This problem aims at modeling situations where a facility provides a service and demand refers to whether a customer requires to be served. Companies providing repair or maintenance services are potential agents in such a modeling framework, and customers may represent sets of users grouped according some criterion (e.g. their location), which are assigned to open facilities that should handle the existing demand. The term ``facility'' should be looked at in a very general way. For instance, we may be referring to a worker or a team.

One potential example concerns mobile health clinics. This type of facility is usually set to assist some specific area or region previously assigned to it. In case the occurring demand is higher than the service capacity, extra personnel is necessary, which may lead to additional costs. Another potential application of the FLPBD can be found in elevator maintenance: each repair or maintenance team must assist a prespecified set of customers in case they call for service.
If the actual demand turns out to be higher than the service capacity (which may indicate maximum time service), then the service still has to be provided; this may call for a temporary relocation of workers from other teams or simply for outsourcing the service to a third party. Other examples of settings fitting the FLPBD include target-oriented advertisement activities, door-to-door product demonstration, etc. In all these cases, potential customers are previously assigned to the facility and may or may not have actual demand, and the total actual demand may turn out to exceed the estimated service capacity.

The FLPBD was motivated and introduced in \citet{Albareda-SambolaFernandezSaldanha:2011} for the particular case when the distributions of the customers' demands are independent, all with the same demand probability. In \citet{Albareda-SambolaFernandezSaldanha:2011}  two different outsourcing policies were considered, and closed forms for the corresponding recourse functions were presented. The obtained numerical results showed that the proposed methodology was computationally highly demanding as the sizes of the instances increased.
In \citet{Albareda-SambolaFernandezSaldanha:2017} the same authors considered the heterogeneous version of the problem where the demand probabilities are not necessarily the same, for the same two outsourcing policies, again under the assumption of independent demands. For that case, evaluating exactly  the cost of a solution becomes computationally unaffordable, so several heuristics, based on GRASP and Path Relinking, were proposed in which solution costs were computed by simulation.
The problem and some of its variants has attracted the attention of other researchers such as \citet{Bieniek:2015}, who considered stochastic i.i.d. demands with arbitrary probability distributions.

In this paper we focus on the analysis of several outsourcing policies for the FLPBD, which guide the two-stage model that we consider.
The first-stage decision is to select a set of facilities (plants) to open together with a \emph{tentative}  (\emph{a priori}) allocation of customers within the set of selected facilities.
The second-stage solution, which determines the recourse action, builds a specific assignment of customers to open plants for each possible realization of the customers' demand. Since each facility has a capacity in terms of the maximum number of customers that it can really serve, for some of the scenarios it may happen that the number of customers who actually have demand, with an \emph{a priori} allocation to a given open plant, exceeds the capacity of the plant. In such cases the outsourcing policy dictates how to re-allocate the exceeding demand, so the total existing demand is eventually served.

In the above framework the choice of the outsourcing policy becomes a highly relevant issue, as it indicates how to \emph{react} in those scenarios when the existing demand exceeds the available capacity, and how to evaluate such a \emph{reaction}. In other words,  the outsourcing policy determines the criterion according to which the quality of solutions will be assessed.
Hence, some outsourcing policies may be better suited than others for some specific classes of scenarios, so different outsourcing policies may lead to different solutions.
Given that one of the main goals of stochastic models is to produce solutions that are \emph{robust} for the considered scenarios, determining a suitable outsourcing policy is on itself a strategic decision, which must be made in advance, and may have a notable impact on the specific location/allocation decisions.
Still, to the best of our knowledge, there is no comparative analysis of alternative recourse functions for stochastic discrete location models, with the exception of \cite{Pages-et-al:2019}. In \cite{Pages-et-al:2019} two recourse functions are considered and compared, exclusively through the cost of their respective optimal solutions.
In this paper we develop an extensive comparative analysis of the performance of several outsourcing policies. On the one hand, we consider four different recourse functions and, on the other hand, the comparison is carried out through the \emph{a priori} decisions produced by each outsourcing policy, which are also evaluated  from the perspective of the other outsourcing policies.

The two outsourcing policies considered in \citet{Albareda-SambolaFernandezSaldanha:2011,Albareda-SambolaFernandezSaldanha:2017} are \emph{facility outsourcing} (FO) and \emph{customer outsourcing} (CO).
Broadly speaking, in FO each facility with insufficient capacity  for serving all its allocated customers who have demand in a given scenario, outsources its deficit of capacity and serves from that plant all its allocated customers with demand.
On the contrary, in CO when in a given scenario some open facility has not enough capacity to serve all the customers allocated  to the plant in the \emph{a priori} solution, only a subset of them are in fact served from that facility, while service to the remaining customers is outsourced, so they are directly served from an external source.
In \citet{Albareda-SambolaFernandezSaldanha:2011,Albareda-SambolaFernandezSaldanha:2017} the outsourced customers are selected with an order driven strategy (OD-CO),  according to the FIFO policy relative to the arrival order of service requests.
Given the widely recognized relevance of outsoucing in the context of production and distribution systems (see, e.g., \citealt{BenarochWebsterKaraz:2012} and \citealt{DolguiProth:2013}), in this paper we consider an additional alternative policy for the selection of outsourced customers within CO, which is based on a cost-minimization criterion, and will be referred to as \emph{cost-driven} (CD-CO).
We also analyze an additional outsourcing policy, \emph{reassignment outsourcing} (RO), derived from the possibility of reassigning some customers allocated to an open facility with insufficient capacity to other open facilities with exceeding capacity.
While FO serves all demand customers from the plant they are allocated to, both CO strategies serve outsourced customers from an external source (third party), whereas RO serves each outsourced customer from one open facility, different from the one the customer was allocated to, but with available capacity.

For each of the considered policies, we present a mathematical programming formulation, which allows us to optimally solve the problem when uncertainty is expressed via a set of scenarios. As we empirically show, some  outsourcing strategies lead to models much more difficult to tackle than others. We carry out extensive computational experiments, whose results we summarize and analyze. Furthermore, the formulations are used to analyze the actual performance of the underlying outsourcing policies. 	
This is accomplished considering both correlated and independent demands so that insights can be gather with respect to possible data dependency.
In particular, we carry out a comparative analysis of the four outsourcing policies, by evaluating the quality of the solutions produced by each specific strategy relative to the other outsourcing policies, i.e., we consider the possibility of using the optimal solution for some outsourcing strategy as an approximate solution for the others.
The results show that for some outsourcing strategies, a good approximation can be obtained by adopting the solution induced by other strategies, and allow us to derive some managerial insight.

In this paper we contribute to existing work in several ways. First, we extend the set of outsourcing policies  considered in \citet{Albareda-SambolaFernandezSaldanha:2011,Albareda-SambolaFernandezSaldanha:2017} with two additional strategies, and propose a mixed-integer linear programming formulation for each of them. Second, we carry out an empirical comparison of the considered outsourcing policies in terms of both the computational performance of the proposed formulations as well as their capability of producing good quality solutions for the other policies. The comparison of the quality of the solutions obtained with the different policies allows us to derive some managerial insight, which could help the decision maker in determining the most suitable outsourcing policy in terms of robustness or other possible indicators.
Finally, our analysis is developed in a general setting where uncertainty is expressed by means of a set of scenarios and it is no longer assumed that the probability distributions of the customers' demands are independent.
Again, this goes beyond the existing literature. The assumption that customers demands are independent holds when customers do not obey to some common interest and demand is not seasonal, but such an assumption is difficult to justify in other cases. Moreover, when ``customers'' and ``facilities'' do not necessarily represent physical entities. This would be the case, for instance, of a situation where ``customers'' were students and ``facilities'' offered courses. Then it would be unlikely that demand be uncorrelated as it could depend, for instance, on the students background.

The paper is organized as follows. In Section \ref{sect:theproblem} we introduce some notation and formally define the FLPBD. Section \ref{policies} focuses on the four alternative outsourcing policies that we have considered. Specifically, we present a mixed-integer programming formulation for each of them, which will be used in the computational experiments. Section \ref{sec:comput} is dedicated to the computational experiments.  The sets of test instances that we have used and their characteristics are described in Section  \ref{subsect:test-instances}. The performance of each of the presented formulations for the different sets of benchmark instances at varying values of the input parameters is summarized and analyzed in Section \ref{subsect:performance}, whereas Section \ref{subsect:cross-comparison} summarizes the results of an extensive comparison of the different policies and derives managerial insigt. 
The paper ends in Section \ref{sec:conclu} with a summary of our main findings and some final comments. 

\section{Definition of the problem} \label{sect:theproblem}
Let $I$ and $J$, with $n=|J|$, denote the set of indices for the potential locations of facilities and for customers, respectively. We assume that the demands of service of customers follow Bernoulli probability distributions, not necessarily independent, with probabilities $p_j, j\in J$. We assume that such uncertainty is expressed by means of a set of possible realizations (scenarios), and denote by $\Omega$ the set of all scenarios, by $\pi^\omega$ the probability of scenario $\omega$ ($\sum\nolimits_{\omega\in \Omega} \pi^\omega =1$), and by  $d_j^\omega \in \{0, 1\}$ the demand of customer $j\in J$  in scenario $\omega\in \Omega$.
Following the terminology introduced in \citet{Albareda-SambolaFernandezSaldanha:2011} the customers with demand in a scenario are referred to as \textit{demand customers} and $D^\omega=\sum_{j\in J}d_j^\omega$ indicates the number of demand customers in scenario $\omega$.

We have the following additional data. For each potential location $i \in I$, $f_i$  is the fixed setup cost for opening facility $i$; $\ell_i$  is a lower bound on the number of customers that have to be {\em assigned} to facility $i$ if it is opened; and, $K_i$  the maximum number of customers that can be {\em served} from facility $i$ if it is opened. For each pair $i\in I, j\in J$, $c_{ij}$ is the cost for serving customer $j$ from facility $i$.

For a given scenario $\omega\in \Omega$  not all the customers need to have demand. This is why we distinguish between the {\em assignment} of customers to plants, which is done {\em a priori} and is independent of the potential realizations, and the {\em service} of customers from open plants, which is decided {\em a posteriori}, once the realization is known.
An {\em a priori} solution is given by a set of {\em operating} (open) facilities together with an assignment of all the customers to these facilities, such that for any open plant the number of customers that are assigned to it is at least $\ell_i$. Since $K_i$ is an upper bound on the number of customers that can be served from an open plant, it does not affect the feasibility of {\em a priori} solutions. Let $i(j)\in I$ denote the facility to which customer $j\in J$ is assigned in the {\em a priori} solution and  $J_i=\{j\in J: i(j) = i\}$, the set of customers assigned to facility $i$ in the {\em a priori} solution.

Given an {\em a priori} solution, the {\em a posteriori} solution indicates the decisions to make once demand customers are known, i.e., it describes the actual services to demand customers.
Let $J_i^{\omega}=J_i\cap\{j\in J: d_j^\omega=1\}$ denote the set of customers assigned to facility $i\in I$ with demand in scenario $\omega$, and $\eta_i^\omega=|J_i^{\omega}|$ the number of such customers.
If the number of demand customers assigned to an open facility $i\in I$ does not exceed its upper bound, i.e. $\eta_i^\omega \le K_i$, then in the {\em a posteriori} solution  all customers indexed in $J_i^\omega$ receive service from plant $i$, each of them incurring a service cost $c_{ij}$, $j\in J$. Instead, when $\eta_i^{\omega} > K_i$ the {\em a posteriori} solution consists of serving $K_i$ (out of $\eta_i^{\omega}$) demand customers from facility $i$ and outsourcing the remaining $\eta_i^{\omega}-K_i$.
A penalty is incurred for every outsourced demand customer. The way in which, for a realization,  it is decided whether a demand customer assigned to a plant with $\eta_i^{\omega}>K_i$, is actually served from $i$ or outsourced, and its corresponding penalty, depends on the outsourcing policy that is applied (see Section \ref{policies} below). The recourse function is the expected cost of the {\em a posteriori} solution, over all possible realizations of the demand vector.

A penalty cost $g_i$ is incurred for every outsourced demand customer. The way in which, for a realization,  it is decided whether a demand customer assigned to a plant with $\eta_i^{\omega}>K_i$, is actually served from $i$ or outsourced, and its corresponding penalty, depends on the outsourcing policy that is applied (see Section \ref{policies} below). The recourse function is the expected cost of the {\em a posteriori} solution, over all possible realizations of the demand vector.

The FLPBD consists of finding a set of facilities to open and an allocation of the customers to the opened facilities, such that the lower bounds $\ell_i$ are satisfied, and the sum of the fixed cost associated with the open
facilities and the recourse function is minimized.

To formulate the FLPBD we define the following sets of decision variables:\\

$y_i = \left\{
\begin{array}{ll}
	1 & \mbox{if a facility is established at } i, \\
	0 & \mbox{otherwise,}
\end{array}
\right. (i \in I).$
\vspace{0.5cm}

$x_{ij} = \left\{
\begin{array}{ll}
	1 & \mbox{if customer $j$ is allocated to $i$,} \\
	0 & \mbox{otherwise,}
\end{array}
\right. (i \in I, j \in J).$\\

The generic formulation for the FLPBD proposed in \citet{Albareda-SambolaFernandezSaldanha:2011} is:
\begin{alignat}{3}
	(P) \quad \min & \quad & & \sum_{i \in I} f_i y_i + {\mathcal Q}(x),& & \label{P:of} \\
	\text{s. t.} \quad & & &  \sum_{i \in I }x_{ij}= 1,              \qquad  & & j \in J,  \label{P:const1}\\
	&&&  x_{ij} \leq y_i,                       \qquad & & i \in I, \: j \in J, \label{P:const2}\\
	&&&  \ell_i y_i \leq \sum_{j \in J} x_{ij}, \qquad & & i \in I,   \label{P:const3}\\
	&&&  y_i \in \{0,1\},                       \qquad & & i \in I,  \label{P:const4}\\
	&&&  x_{ij} \in \{0,1\},                    \qquad & & i \in I, \: j \in J.  \label{P:vars}
\end{alignat}
The objective function \eqref{P:of} includes the fixed costs for opening the facilities and the recourse function. In particular,
${\mathcal Q}(x) = \mathbb{E} \left[ \mbox{Service cost}+ \mbox{Penalty cost} \right]$. Constraints \eqref{P:const1} assure that all customers will be assigned to (exactly) one facility while constraints \eqref{P:const2} impose that these assignments are only done
to operating facilities.
Constraints \eqref{P:const3} state the minimum number of customers that must be assigned to each operating facility. Finally, \eqref{P:const4}--\eqref{P:vars} define the domain of the variables.

\section{Outsourcing Policies}\label{policies}
The expression of $\mathcal {Q}(x)$ and the additional variables and constraints that may be needed to express the FLPBD through a mathematical programming formulation are directly related to how the recourse action is defined.
That is, what specific outsourcing policy is applied. Below we describe the outsourcing policies that will be considered in this work and we discuss manageable Mixed-Integer Linear Programming (MILP) formulations in each case. Broadly speaking, the alternative policies that we study differ from each other in the recourse actions that are taken in the scenarios where some facility has a number of assigned demand customers that exceeds its capacity.

\subsection{Facility outsourcing}
With the facility outsourcing (FO) policy, under scenario $\omega$, facility $i$ takes delivery of the whole set $J_i^{\omega}$. When $\eta_i^{\omega}>K_i$, then $\eta_i^{\omega}-K_i$ units of product are
outsourced, at a unit cost $g_i$. Then, facility $i$ serves the full demand of its assigned customers, $\eta_i^{\omega}$, at the same cost $c_{ij}$ that would be incurred if it were not outsourced.
This recourse action models an external purchase of the resources needed to fully satisfy the demand of an open facility and was applied in \citet{Albareda-SambolaFernandezSaldanha:2011} for the case when the demand of customers are independent random variables, and all of them have the same probability of demand, i.e., $p_j=p$.

With the FO policy, each scenario $\omega\in \Omega$ is characterized by its probability $\pi^{\omega}$ and the demands $d_j^\omega \in \{0, 1\}$, $j\in J$.
To formulate the FO-FLPBD, in addition to the $y$ and $x$ decision variables introduced above, we use the following:

\begin{itemize}
	\item[] $\theta_i^\omega:$ number of demand customers outsourced at facility $i$ under scenario $\omega$.  ($i \in I, \omega\in \Omega$).\\
	
	\item[] $z^\omega:$ total penalty incurred under scenario $\omega \in \Omega$.\\
\end{itemize}

The formulation is:
\begin{alignat}{2}
	\mbox{FO} \hspace*{0,3cm} \min \quad & \sum_{i \in I} f_i  y_i + \sum\limits_{i\in I} \sum\limits_{j\in J} p_j c_{ij} x_{ij} +
	\sum\limits_{\omega \in \Omega} \pi^\omega z^\omega, && \label{P1:of} \\
	\text{s. t.}  \quad &\eqref{P:const1}-\eqref{P:const3}, &&\nonumber\\
	&   \theta_i^\omega \geqslant \sum\limits_{j\in J} d_j^\omega x_{ij} - K_i{y_i},        &\qquad& i \in I, \omega \in \Omega, \label{P1:noutsourced}\\
	& z^\omega \geqslant \sum\limits_{i\in I} g_i\theta_i^\omega, \quad         && \omega \in \Omega, \label{P1:speed1}\\
	&  z^\omega\geqslant 0,  \theta_{i}^\omega \geqslant 0,                                   && i \in I, \omega\in \Omega, \label{P1:thetapos}\\
	&  y_i \in \{0,1\},                                                                       && i \in I, \label{P1:ybin}\\
	&  x_{ij} \in \{0,1\},                                                                   && i \in I, j \in J.  \label{P1:xbin}
\end{alignat}
The objective function \eqref{P1:of} includes the costs for opening facilities plus the expected value of the service plus outsourcing costs.
As explained, Constraints  \eqref{P:const1}--\eqref{P:const3} guarantee the feasibility of the \emph{a priori solution}.
Constraints \eqref{P1:noutsourced} force $\theta$ variables to take consistent values, and Constraints \eqref{P1:speed1} compute the penalty cost of each scenario.
Indeed \eqref{P1:speed1} will hold as equality in any optimal solution.
So, in fact, they are not strictly needed, as their right hand side could be substituted in the last term of the objective function instead of using the variables $z^\omega$.
Some preliminary experiments showed that the formulation with Constraints \eqref{P1:speed1} could be solved in smaller times than that where its right was substituted in the objective function. Hence, we used this alternative, even if we have no theoretical argument that justifies this improvement.
The domain of the variables is defined by \eqref{P1:thetapos}--\eqref{P1:xbin}.
Note that, given the structure of the formulation, integrality constraints on the $z$ and $\theta$ variables can be relaxed to nonnegativity constraints.

Formulation \eqref{P1:of}--\eqref{P1:xbin} uses $|I|(1+|J|)+|\Omega|(|I|+1)$ variables and has $|J|(1+|I|)+|I|+|\Omega|(|I|+1)$ constraints.
Depending on the size of $\Omega$, these numbers can be quite high, even for moderate numbers in terms of customers and facilities.
Hence, enhancing the formulation can be very useful to
decrease the CPU time required to solve such model to proven optimality using an off-the-shelf solver.
Inequalities \eqref{P1:speed2} and \eqref{P1:speed3} below have proven to give a good balance between the increase in the size of the formulation and the improvement obtained when solving the model.
\begin{alignat}{2}
	&   \sum\limits_{\omega\in \Omega} \pi^\omega \theta_i^\omega \geqslant \sum\limits_{j\in J}p_j x_{ij} - K_i y_i, &\qquad& i\in I, \label{P1:speed2}\\
	&   \sum\limits_{i\in I}K_iy_i + \sum\limits_{i\in I} \theta_i^{\tilde{\omega}} \geqslant D^{\tilde{\omega}}. &&  \label{P1:speed3}
\end{alignat}
Inequality \eqref{P1:speed2} states that the expected number of demand customers outsourced at facility $i$ ($i \in I$) is at least the expected number of demand customers assigned to that facility minus the capacity of the facility. Note that these inequalities can be derived as a weighted sum over all scenarios of Constraints \eqref{P1:noutsourced},  using as weights $\pi^\omega$, $\omega\in \Omega$, after imposing that $\sum\limits_{\omega\in \Omega} \pi^\omega=1$. Since this probability equality affects the input data only but is not explicit in the formulation, Constraints \eqref{P1:speed2}  are valid inequalities, which are not implied by \eqref{P1:noutsourced}, and may help cutting fractional solutions.
Note that, for binary solutions, these constraints are activated only if $y_i=1$.

In \eqref{P1:speed3}, ${\tilde{\omega}}\in \Omega$ is the scenario with the largest number of demand customers ($D^{\tilde{\omega}}$).
This constraint ensures that the maximum number of customers that can be served by the open facilities plus the outsourced demand customers, is never below the total demand.
This constraint holds for every scenario, but adding such a constraint for all scenarios would increase considerably the size of the formulation, which explains why we consider solely \eqref{P1:speed3}. Even being one single constraint, it has proven to decrease the computation time in some instances.

\subsection{Customer outsourcing}

With the customer outsourcing (CO) strategy, in the scenarios where the number of demand customers assigned to facility $i$ exceeds its capacity $K_i$, i.e.,  $\eta_i^\omega>K_i$, exactly $K_i$ customers are served from facility $i$, whereas the remaining $\eta_i^\omega-K_i$ customers of $J_i^\omega$ are outsourced and receive service from an external third party.
Service costs $c_{ij}$ are incurred for the customers served from facility $i$, whereas a penalty $g_i$ is incurred for each outsourced customer, which depends on the facility $i$ the customer is assigned to.
Hence, to formulate the FLPBD with a CO policy new decision variables are needed, in addition to the ones defined above, to identify the outsourced customers. In particular, we define:\\

$s_{ij}^\omega= \left\{
\begin{array}{ll}
	1 & \mbox{if customer } j \textrm{ is served from facility } i \textrm{ under scenario } \omega, \\
	0 & \mbox{otherwise},
\end{array}
\right.$ 

\medskip
\noindent $i \in I, j\in J, \omega \in \Omega$.\\

We consider two different versions of the CO policy. In the first one, that we call cost-driven CO policy (CD-CO), the decision of whether a demand customer $j\in J_i^{\omega}$ is served from  facility $i$ or outsourced is based solely on a cost-minimization criterion.
Thus, similarly to the FO policy, in the CD-CO policy a scenario $\omega\in \Omega$ is characterized by its probability $\pi^{\omega}$ and the demands $d_j^\omega \in \{0, 1\}$, $j\in J$.
The formulation for the FLPBD  with CD-CO is:
\begin{alignat}{2}
	\mbox{CD-CO} \hspace*{0,3cm} \min &  \sum_{i \in I} f_i  y_i + \sum\limits_{\omega \in \Omega} \sum\limits_{i\in I} \sum\limits_{j\in J} c_{ij}s_{ij}^\omega +&& \sum\limits_{\omega \in \Omega} \pi^\omega z^\omega,            \label{P2:of} \\
	s.\:t.  \qquad &\eqref{P:const1}-\eqref{P:const3}, &&\nonumber\\
	&\sum\limits_{\omega \in \Omega}   s_{ij}^\omega \leqslant \left(\sum\limits_{\omega\in \Omega} d_j^\omega\right) x_{ij}, && i \in I, j \in J,  \label{P2:link_sx}\\
	&   \sum\limits_{j\in J} d_j^\omega s_{ij}^\omega \leqslant K_i,                         && i \in I, \omega \in \Omega, \label{P2:capacity}\\
	&   \sum\limits_{j\in J} d_j^\omega s_{ij}^\omega + \theta_i^\omega {\geqslant} \sum\limits_{j\in J} d_j^\omega x_{ij}, && i \in I, \omega \in \Omega, \label{P2:noutsourced}\\
	& z^\omega \geqslant \sum\limits_{i\in I} g_i \theta_i^\omega,\quad         && \omega \in \Omega, \label{P2:speed1}\\
	&  z^\omega\geqslant 0,  \theta_{i}^\omega \geqslant 0,                        && i \in I, \omega\in \Omega, \label{P2:thetapos}\\
	&  y_i \in \{0,1\},                                                                      && i \in I,\label{P2:ybin}\\
	&  x_{ij} \in \{0,1\},                                                                   && i \in I, j \in J,  \label{P2:xbin}\\
	&  s_{ij}^\omega \in \{0,1\},                                                            && i \in I, j \in J, \omega\in \Omega.  \label{P2:sbin}
\end{alignat}

Again, Constraints  \eqref{P:const1}--\eqref{P:const3} guarantee the feasibility of the \emph{a priori} solution.
The second stage variables $s_{ij}^\omega$ are now used to compute the expected service cost. Constraints \eqref{P2:link_sx} ensure that service from open facilities is only provided according to the \emph{a priori} assignments dictated by the $x$ variables.
Constraints \eqref{P2:capacity} and \eqref{P2:noutsourced} state the service capacities of the facilities and set the right value of the number of outsourced units at each facility, respectively.
Again, the structure of the problem allows to relax integrality constraints on $\theta$ variables and restrict them to be just nonnegative.

The number of variables of formulation CD-CO  has increased in $|I|\times|J|\times|\Omega|$ with respect to the number of variables of the FO formulation. Its number of constraints is also larger, as it has raised to $|J|(1+2|I|)+|I|+|\Omega|(2|I|+1)$.

In the second CO strategy that we consider, in the scenarios where $\eta_i^\omega>K_i$, the demand customers to serve from facility $i$ are selected following a FIFO policy, relative to the order in which requests of service have arrived. This order driven recourse action will be referred to as OD-CO and was applied in \citet{Albareda-SambolaFernandezSaldanha:2011} for the particular case when $(i)$ the demand of customers are independent random variables, $(ii)$ all of them have the same probability of demand, i.e., $p_j=p$, and $(iii)$ requests of service arrive in a random order.

Note that for the FLPBD with OD-CO policy a scenario $\omega\in \Omega$  is no longer fully characterized by its probability and the demands of the customers. Now, the order in which calls for service from the demand customers of $J_i^\omega$ were received must also be known.

The formulation for the OD-CO is \eqref{P2:of}--\eqref{P2:sbin} plus the following set of constraints, which ensure that the FIFO policy is followed for selecting the customers that will be served from a given facility when $\eta^{\omega}_i>K_i$.
The notation $j' \prec^\omega j$ indicates that customers $j, j' \in J$ have demand under scenario $\omega$  and $j'$ requested service before $j$.
\begin{alignat}{2}
	& K_i d_j^\omega( x_{ij} - s_{ij}^\omega ) \leqslant \sum\limits_{j' \prec^{\omega} j} d_{j'}^\omega s_{ij'}^\omega, \qquad && i \in I, \omega \in \Omega, j\in J.\label{P2:order}
\end{alignat}

\subsection{Reassignment outsourcing}
In the previous policies, when facility $i$ has a number of assigned demand customers that exceeds its capacity $K_i$, then the excess of demand at $i$ is served from external sources.
Instead, in the reassignment outsourcing (RO) policy the unused capacity of other open plants can be used to satisfy the deficit of capacity at $i$ prior to resorting to outsourcing. Similarly to the CO policy, in the RO recourse, when $\eta_i^\omega>K_i$ exactly $K_i$ customers of $J_i^\omega$ are served from $i$. Now, up to $\sum_{i'\in I:i'\ne i}\left(K_{i'}-\eta_{i'}^\omega\right)^+$ customers of $J_i^{\omega}$ can be served from different open facilities (with $(a)^+$ standing for $\max\{0,a\}$).
The service cost for serving a demand customer $j\in J_i^\omega$ from an  open facility $i'\in I$ is $c_{i'j}$, independently of whether or not $i(j)=i'$. When $i(j)\ne i'$, i.e. $j$ is not assigned to $i'$ in the \emph{a priori} solution, an additional penalty $h_j$ is paid. Finally, when $\eta_i^\omega-K_i>\sum_{i'\in I:i'\ne i}\left(K_{i'}-\eta_{i'}^\omega\right)^+$ we resort to external outsourcing to serve the remaining unserved demand assuming a unit cost $g$ per outsourced customer.

RO is a cost-driven policy, so the decision of whether a demand customer is served from its allocated facility, a different open facility, or outsourced is only based on the cost-minimization criterion.

In the RO-FLPBD, each scenario $\omega\in \Omega$ is characterized by its probability $\pi^\omega$ and the customers demands $d_j^\omega \in \{0, 1\}$, $j\in J$. In this case, to identify the demand customers served from facilities different from the ones they are assigned to and the demand customers who receive service from an external source, we use the following sets of binary decision variables in addition to the location ($y$) and the allocation ($x$) variables:

\begin{itemize}
	\item $\lambda_j^\omega =1 \Leftrightarrow$ demand customer $j\in J^\omega$ is reassigned to an open facility $i'\ne i(j)$ in scenario $\omega$.
	
	\item $\mu_j^\omega  =1 \Leftrightarrow$ demand customer $j\in J^\omega$ is served from external sourcing in scenario $\omega$.
\end{itemize}

The obtained formulation is:
\begin{alignat}{3}
	\mbox{RO} \quad \min & \quad & &  \sum_{i \in I} f_i  y_i + \sum\limits_{\omega \in \Omega} \pi^\omega \sum\limits_{j\in J} \left( \left(  \sum_{i \in I} c_{ij} s_{ij}^\omega \right)  + h_j\lambda_j^\omega + g\mu_j^\omega \right),   & \quad & \label{PP1:of}
\end{alignat}
\vspace{-7mm}
\begin{alignat}{3}	
	\rm{s.\:t.}  \qquad & & &  \eqref{P:const1}-\eqref{P:const3}, &&\nonumber\\
	&& &    \sum\limits_{j\in J}  s_{ij}^\omega \leqslant  K_i y_i,                              &\qquad& i \in I, \omega \in \Omega, \label{PP1:noutsourced}\\
	&& & d_j^\omega\lambda_j^\omega \geqslant s_{ij}^\omega - d_j^\omega x_{ij}^\omega,                              && i\in I, j\in J, \omega \in \Omega, \label{PP1:activate1}\\
	&&& \sum\limits_{i\in I} s_{ij}^\omega + \mu_j^\omega \geqslant d_j^\omega,                     && j\in J, \omega \in \Omega, \label{PP1:activate2}\\
	& && y_i , x_{ij}, s_{ij}^\omega, \lambda_j^\omega, \mu_j^\omega \in \{0,1\},          && i \in I, j \in J, \omega\in \Omega. \label{PP1:thetapos}
\end{alignat}

The role of Constraints \eqref{P:const1}--\eqref{P:const3} has been repeatedly explained. Constraints  \eqref{PP1:noutsourced} impose that no facility serves more than $K_i$ customers. Observe that as opposed to the previous policies, these constraints do not prevent that a facility serves a customer not assigned to it in the {\em a priori} solution. Constraints  \eqref{PP1:activate1} activate the $\lambda$ variables associated with demand customers served from open facilities different from the ones they are assigned to in the {\em a priori} solution.
The $\mu$ variables associated with customers that receive service from an external source are activated in constraints  \eqref{PP1:activate2}.

\section{Computational Experiments}\label{sec:comput}

\subsection{Test instances}\label{subsect:test-instances}

All the instances used in this paper are generated from the homogeneous instances presented in \citet{Albareda-SambolaFernandezSaldanha:2011}. In that work, homogeneous FLPBD instances were generated
taking as a starting point
11 Traveling Salesman Problem (TSP) instances available at \url{http://comopt.ifi.uni-heidelberg.de/software/TSPLIB95/} namely, \verb|berlin52|, \verb|eil51|, \verb|eil76|, \verb|kroA100|, \verb|kroB100|, \verb|kroC100|, \verb|kroD100|, \verb|kroE100|, \verb|pr76|, \verb|rat99|, and \verb|st70|.
From those TSP instances, small and large FLPBD instances were generated in \citet{Albareda-SambolaFernandezSaldanha:2011}, with $|I|=15$, $|J|=30$ and $|I|=20$, $|J|=60$, respectively.
Below we briefly recall the generation process.

From each original TSP instance with $N$ nodes, first the number of facilities and customers was set in such a way that $|I|+|J| \leq N$.
This means that the largest instances could be generated only from the TSP instances \verb|kroA100|, \verb|kroB100|, \verb|kroC100|, \verb|kroD100|, \verb|kroE100|, and \verb|rat99|.\\
Then, for each fixed FLPBD dimensions, three different sets of customers and potential plants were randomly selected among the $n$ original nodes.
For each choice of plants and customers, the remaining data for the FLPBD was generated varying several parameters as follows: (i) three different values for the probability of demand ($0.1$, $0.5$ and $0.9$), (ii) three different levels of variability for the setup costs ($0$, $\mu/10$, and $\mu/3$, where $\mu$ is the expected value set for the setup costs), (iii) low, medium and high capacity $K_i$, and (iv) two different possibilities for
$\ell$---the minimum number of customers to assign \textit{a priori} to the opened facilities ($0$ and a value greater than $0$).

In total, \citet{Albareda-SambolaFernandezSaldanha:2011} generated 2754 instances divided into 17 groups (11 groups of small instances and 6 groups of large instances).
A subset of those instances was then considered in \citet{Albareda-SambolaFernandezSaldanha:2017} for generating a data set for the non-homogeneous case.
Three types of customers were considered: low-, medium-, and high-probability demand customers, with demand probabilities  drawn from U(0.10, 0.25), U(0.40, 0.60), and U(0.75, 0.90) distributions, respectively.
In its turn,  different patterns were defined for the overall demand.  In this work we use the following two patterns: in pattern 1, there are 20\% of low-probability demand customers, 60\% medium, and 20\% high and in pattern 2, these percentages are 20\%, 20\%, and 60\%, respectively.

The above procedures generated a whole set of probabilistically defined FLPBD instances (all the details can be found in \citet{Albareda-SambolaFernandezSaldanha:2017}).
For the current work, we have used a subset of the instances with patterns 1 and 2 referred to as \emph{PT1} and \emph{PT2}, respectively.
In particular, we considered the set of instances corresponding to: (i) the intermediate value for the variability for the setup costs ($\mu/10$), (ii) low and high capacity levels, and (iii) $\ell>0$.
Concerning the capacities we recall that the data for every instance in \citet{Albareda-SambolaFernandezSaldanha:2011} was generated considering among other parameters, one multiplication factor---$\gamma \in \{1,2,4\}$---that leads to the so-called low, intermediate and high capacity levels for the facilities.
In this work, we consider the instances generated using $\gamma=1,4$ which means that we are retrieving the instances with the lower and the higher capacities generated in the above mentioned work.

The above choices for the testbed instances to be used in the current work are motivated by the fact that we wanted to consider a rather limited set of instances. In fact, considering all the instances worked out by \citet{Albareda-SambolaFernandezSaldanha:2017} was not relevant for a paper whose focus is not in algorithmic developments. Hence, we concentrated our experiments on instances corresponding to one possibility for the variability of the setup costs (thus we chose the intermediate one), and two possibilities for the capacities: low and high. This led to a base set of 204 instances.
From each of them, we generated two scenario-defined FLPBD instances with $50$ scenarios each, yielding two different groups of 204 instances each. The number of scenarios was fixed to get a good compromise between quality of randomness representation, and affordability in terms of computational effort.
In the first group, for each scenario, demands were generated independently, and according to the previously defined probabilities.
After that, the customers having demand are randomly sorted to simulate their calling sequence. In the second one, we used again the previous demand probabilities, but  now forcing spatially correlated demands. To do so, for each pair of customers $j,j'\in J$ we first computed $\delta_{jj'}=\max\{0.1, \max_{i\in I}\{|c_{ij}-c_{ij'}|\}\}$ which, although not being a proper distance, gives an idea of the proximity of two customers, taking the set of facilities as a reference. The idea is to force that the demands of customers that are at a small {\em distance} have higher correlations. The decay of this correlation is governed by a preset decreasing function $w(\delta)$. In particular, we used $w(\delta)=1+\left(1-2\frac{\delta}{\Delta}\right)^\frac{1}{3}$, where $\Delta$ stands for the maximum distance between two customers. Using the above {\em distances} and decay function, we generated each scenario using the following observation:  Let $\{Z_j\}_{j\in J}$  be a family of independent random variables, each following a $N(0; 1)$ distribution. Then variables $\{Y_j\}_{j\in J}$ defined as
$$ Y_j = \frac{ \sum\limits_{k\in J} w(\delta_{jk})Z_k}{\sqrt{\sum\limits_{k \in J} w(\delta_{jk})^2}},$$
also follow a $N(0; 1)$ distribution and have positive correlations which increase as the corresponding customers get closer.
Then, from a sample of the $\{Z_j\}_{j\in J}$ variables we computed the corresponding $\{Y_j\}_{j\in J}$ values and set
$$
d_j=\begin{cases}1 & \textrm{if } \Phi(Y_j)\leqslant p_j, \\
	0 & \textrm{otherwise},\end{cases} \textrm{\quad  with } \Phi \textrm{ denoting the cdf of the } N(0; 1) \textrm{ distribution.}
$$
By proceeding like this, we generated demands following correlated Bernoulli distributions with the same probabilities as in the instances with independent demands. Again, after generating these demands we simulate the calling sequence by randomly sorting the customers with demand.

In total we have 408 instances (204 with uncorrelated demands and 204 with correlated demands).
In the following, ``U'' stands for instances with uncorrelated demands and ``C'' for instances with correlated demands.
Recall also that $\gamma=1$ indicates low capacities in use whereas $\gamma=4$ corresponds to high capacites.
Table~\ref{tab:instances} summarizes the characteristics of the scenario-defined instances
($|\Omega|=50$).
\begin{table}[!h]
	\centering
	\begin{tabular}{c|c|cc|r}
		&             &  \texttt{PT1}       &   \texttt{PT2}      &   Total \\
		& $\gamma$            & (20\%, 60\%, 20\%)  &  (20\%, 20\%, 60\%)  &\\
		\hline
		Small              & 1    &    33 (U) + 33 (C)     &        33 (U) + 33 (C)         & 132 \\
		$|I|=15, |J|=30$   & 4    &    33 (U) + 33 (C)     &        33 (U) + 33 (C)         & 132 \\
		\hline
		Large              & 1    &   18 (U) + 18 (C)      &       18 (U) + 18 (C)          &  72\\
		$|I|=20, |J|=60$   & 4    &   18 (U) + 18 (C)      &       18 (U) + 18 (C)          &  72\\
		\hline
		Total       &               &       204       &         204             & 408 \\
	\end{tabular}
	\caption{Summary of scenario-defined instances. \label{tab:instances}}
\end{table}

All the formulations presented in this paper were implemented using CPLEX 12.6 callable libraries within a C code. In all cases a CPU time limit was set. This limit was one hour for the small instances, and two hours for the large ones. All the tests were carried out on a Pentium(R) 4, 3.2GHz, 1.0GB of RAM.

In what follows, we report the results obtained. First we analyze the performance of the formulations for the different outsourcing policies and then we carry out a comparative analysis among them.

\subsection{Performance of formulations}\label{subsect:performance}

In Tables~\ref{table:summary:fac_out}--\ref{table:RO} we summarize the results obtained when using the four models proposed in the previous section.
In these tables, ``S'' and ``L'' stand for small and large instances, respectively. The columns under \emph{CPU} give average computing times in seconds. The columns under $\%gap$ give average percent optimality gaps, computed, for each instance, as $100\frac{z_U-z_L}{z_L}$, where $z_U$ and $z_L$ respectively denote the values of the best feasible solution and the best lower bound at termination. In each case, averages are computed over all the instances of the corresponding row. Columns under \emph{\# opt} indicate the number of proven optimal solutions found. The last row in each table summarizes the results over the full set of benchmark instances, giving average computing times in columns CPU and $\%gap$, respectively, and total number of optimal solutions found in \emph{\# opt}. The other entries have already been explained before.

\begin{table}[!h]
	\centering
	\begin{tabular}{|c|c|c|rr|rr|ccc|}
		\hline
		&&	  &              \multicolumn{2}{c|}{CPU} &  \multicolumn{2}{c|}{\%gap}   & \multicolumn{3}{c|}{\# opt}\\		
		&&	 $\gamma$ &                PT1	     &       PT2	  &     PT1	 &   PT2  &   PT1	 &    PT2	 &   Total \\
		\hline
		\multirow{4}{*}{U}  &    \multirow{2}{*}{S} &  1    &      520.63  &        434.19    &     0.05    &     0.03    &      31    &     31    &     62     \\
		&                       &  4    &        5.89  &          3.70    &     0.00    &     0.00    &      33    &     33    &     66       \\
		\cmidrule{2-10}
		&    \multirow{2}{*}{L} &  1    &     2467.54  &       2328.09    &     2.76    &     1.11    &       1    &      0    &      1  \\
		&                       &  4    &       81.78  &         62.23    &     0.00    &     0.00    &      18    &     18    &     36    \\
		\hline
		\multirow{4}{*}{C}  &    \multirow{2}{*}{S} &  1    &        0.33  &          0.37    &     0.00    &     0.00    &      33    &     33    &     66 \\
		&                       &  4    &        1.36  &          0.45    &     0.00    &     0.00    &      33    &     33    &     66 \\
		\cmidrule{2-10}
		&    \multirow{2}{*}{L} & 1    &        2.82  &          2.10    &     0.00    &     0.00    &      18    &     18    &     36 \\
		&                       &  4    &        3.91  &          1.95    &     0.00    &     0.00    &      18    &     18    &     36 \\
		\hline
		Summary               &                       &       &      310.98  &        282.23    &     0.25    &     0.10    &     185    &    184    &    369 \\
		\hline
	\end{tabular}
	\caption{Summary of results for model FO---facility outsourcing policy.}\label{table:summary:fac_out}
\end{table}

Table~\ref{table:summary:fac_out} contains the results for formulation FO.
Looking into this table, we observe that the instances with uncorrelated demands are harder to handle than those with correlated demands.
This is expected since a correlated behaviour somehow simplifies (reduces) the range of observable ``futures'' thus leading to easier-to-solve instances.
The superiority of the correlated instances is clear both in terms of the number of instances solved to optimality and in terms of the computing time required by the solver to accomplish that.
As expected, the larger instances are harder to tackle: they require on average a much higher computing time and also more space, and thus fewer are being solved to proven optimality within the time limit. It is worth noting that, as will be seen in the next section, most of the non-solved instances terminated because of memory problems. Actually, only two small instances terminated because of the time limit.
When using formulation FO, the demand pattern considered does not seem to influence the results obtained.
This indicates that having a higher or lower expected service request is not as relevant as having correlated or uncorrelated demand or higher or lower capacities for the facilities. Concerning the capacity type for the instances with uncorrelated demands, it clearly has an effect: the instances with tighter capacities ($\gamma =1$) are more challenging.
Overall, when using model FO there is empirical evidence that the harder instances are those with uncorrelated demands and low capacities.

In Table~\ref{table:summary:CO_FIFO} we can find a synthesis of the results obtained when considering formulation OD-CO.
Clearly, the structure of the mathematical model derived for this outsourcing strategy makes it more challenging. The effect of the correlation in the demand does not impact as much as before in terms of the CPU time required to solve the instances to proven optimality.
However, for correlated demands, more instances are solved successfully. In turn, the size of the instances and the type of capacity do influence the results. In fact, similarly to formulation FO, tighter capacities ($\gamma=1$) make the instances harder to tackle.
This holds in terms of the number of instances solved to proven optimality as well as in terms of the computing times required to solve such instances and in terms of the final gap for the remaining instances. As observed in the case of facility outsourcing, the demand pattern does not seem to impact on the results. Overall, when considering model OD-CO, the harder instances to tackle are those with uncorrelated demands and tighter capacities.

\begin{table}[!h]
	\centering
	\begin{tabular}{|c|c|c|rr|rr|ccc|}
		\hline
		&&	  &              \multicolumn{2}{c|}{CPU} &  \multicolumn{2}{c|}{\%gap}   & \multicolumn{3}{c|}{\# opt}\\		
		&& $\gamma$ &                PT1	     &       PT2	  &     PT1	 &   PT2	  &   PT1	 &   PT2	 &   Total \\
		\hline
		\multirow{4}{*}{U}  &    \multirow{2}{*}{S} &  1    &     3600.12  &       3600.09    &     4.96    &     4.10    &      0    &     0    &      0     \\
		&                       &  4    &     1704.22  &       1985.06    &     0.47    &     0.65    &     26    &    21    &     47       \\
		\cmidrule{2-10}
		&    \multirow{2}{*}{L} &  1    &     7200.21  &       7200.27    &     7.95    &     8.59    &      0    &     0    &      0  \\
		&                       &  4    &     6786.20  &       6370.27    &     1.34    &     1.14    &      2    &     5    &      7    \\
		\hline
		\multirow{4}{*}{C}  &    \multirow{2}{*}{S} &  1    &     2162.91  &       2584.14    &     1.66    &     3.70    &     10    &     3    &     13 \\
		&                       &  4    &      902.61  &        230.51    &     0.32    &     0.00    &     26    &    33    &     59 \\
		\cmidrule{2-10}
		&    \multirow{2}{*}{L} &  1    &     7145.92  &       6544.58    &     1.23    &     0.63    &      0    &     3    &      3 \\
		&                       &  4    &     4457.38  &       4573.73    &     1.08    &     0.74    &     10    &     9    &     19 \\
		\hline
		Summary               &                       &       &     3611.86  &       3537.22    &     2.22    &     2.35    &     74    &    74    &    148 \\
		\hline
	\end{tabular}
	\caption{Summary of results for model OD-CO---customer outsourcing with a order driven discipline.} \label{table:summary:CO_FIFO}
\end{table}

Table~\ref{table:summary:CO-CD} contains the results for formulation CD-CO.
As for the above models, the size of the instances and the capacity type (low vs high) of the facilities clearly influence the results and in a similar way. However, contrary to the previous cases, now, the demand pattern influences the final gap observed for the instances not solved to optimality.
Although in some cases this is not significant, in the columns headed by ``\%gap'' we observe that averages corresponding to PT2 tend to be smaller than those for PT1.
The demand pattern 2 also seems to make the instances more tractable when it comes to solve them to optimality as we observe in the last columns of the table.
The above aspects may be justified by the fact that in the case of pattern 2 we have more customers with high probability of requesting the service, i.e., we have a larger expected service request. This feature embedded in a cost-driven outsourcing strategy is apparently helping discarding sooner many solutions that are not as competitive and they could be when cost is not a driving feature for deciding about customer outsourcing.
Observing this table we also realize that no instance combining uncorrelated demands with tight capacities ($\gamma=1$) was solved to proven optimality within the time limit. Overall, from the perspective of using model CD-CO, the hardest instances are those with uncorrelated demands, tighter capacities and the first demand pattern (60\% medium-demand customers and 20\% of low- and high-demand customers).

\begin{table}[!h]
	\centering
	\begin{tabular}{|c|c|c|rr|rr|ccc|}
		\hline
		&&	  &              \multicolumn{2}{c|}{CPU} &  \multicolumn{2}{c|}{\%gap}   & \multicolumn{3}{c|}{\# opt}\\		
		&&	 $\gamma$ &                PT1	     &    PT2		  &    PT1		 &   PT2	  &   PT1		 &    PT2	 &   Total \\
		\hline
		\multirow{4}{*}{U}  &    \multirow{2}{*}{S} &  1    &     2783.10  &       3277.74    &     2.56    &     1.95    &       2    &     3    &     5     \\
		&                       &  4    &      580.90  &        560.63    &     0.12    &     0.11    &      30    &    31    &    61       \\
		\cmidrule{2-10}
		&    \multirow{2}{*}{L} &  1    &     6368.24  &       6930.04    &     6.14    &     3.40    &       0    &     0    &     0  \\
		&                       &  4    &     3746.21  &       2394.78    &     0.09    &     0.39    &      14    &    16    &    30    \\
		\hline
		\multirow{4}{*}{C}  &    \multirow{2}{*}{S} &  1    &      852.81  &       1422.87    &     0.20    &     0.47    &      29    &    26    &    55 \\
		&                       &  4    &      358.68  &         22.55    &     0.31    &     0.00    &      27    &    33    &    60 \\
		\cmidrule{2-10}
		&    \multirow{2}{*}{L} &  1    &     2951.89  &       1036.39    &     0.21    &     0.06    &      15    &    17    &    32 \\
		&                       &  4    &     1326.28  &        655.00    &     0.36    &     0.23    &      14    &    16    &    30 \\
		\hline
		Summary               &                       &       &     2010.09  &       1826.75    &     1.12    &     0.77    &     131    &   142    &   273 \\
		\hline
	\end{tabular}
	\caption{Summary of results for model CD-CO---customer outsourcing using a cost-driven discipline.} \label{table:summary:CO-CD}
\end{table}

\begin{table}[!h]
	\centering
	\begin{tabular}{|c|c|c|rr|rr|ccc|}
		\hline
		&&	  &              \multicolumn{2}{c|}{CPU} &  \multicolumn{2}{c|}{\%gap}   & \multicolumn{3}{c|}{\# opt}\\		
		&& $\gamma$ &          PT1	     &     PT2	  &    PT1	 &  PT2	  &   PT1	 &    PT2	 &   Total \\
		\hline
		\multirow{4}{*}{U}  &    \multirow{2}{*}{S} &  1    &     3062.82  &       3006.35    &     1.18    &     0.56    &      0    &      0    &      0     \\
		&          &  4    &     1082.92  &       1173.42    &     0.07    &     0.07    &     23    &     24    &     47       \\
		\cmidrule{2-10}
		&    \multirow{2}{*}{L} &  1    &     7200.12  &       7200.18    &     3.05    &     2.27    &      0    &      0    &      0  \\
		&          &  4    &     7099.35  &       7200.13    &     0.46    &     0.70    &      0    &      0    &      0    \\
		\hline
		\multirow{4}{*}{C}   &    \multirow{2}{*}{S} &  1    &      558.38  &        84.18     &     0.02    &     0.01    &     30    &     32    &     62 \\
		&          &  4    &      226.63  &       126.29     &     0.01    &     0.01    &     32    &     32    &     64 \\
		\cmidrule{2-10}
		&    \multirow{2}{*}{L} &  1    &     5401.62  &       3965.66    &     0.13    &     0.06    &      3    &      7    &     10 \\
		&          &  4    &     4740.20  &       3484.75    &     0.08    &     0.04    &      7    &     11    &     18 \\
		\hline
		Summary   &          &       &     2954.21  &       2638.19    &     0.53    &     0.38    &     95    &    106    &    201 \\
		\hline
	\end{tabular}
	\caption{Summary of results for model RO---reassignment outsourcing.}\label{table:RO}
\end{table}

Finally, in Table~\ref{table:RO} we present the results for the fourth outsourcing strategy proposed: reassignment outsourcing.
The conclusions are very similar to those drawn for model CD-CO.
Again, the demand pattern seems to have an influence in the results and in a similar way as for the previous model discussed.
Now, no large instance with uncorrelated demands could be solved to optimality within the time limit.
In the case of loose capacities and for the first demand pattern we observe an average computing time of 7099 seconds, i.e., below the time limit.
This is an indication that for some instances in that group the machine ran out of memory before the time limit was reached. Actually, if we only consider the instances that could be optimally solved, the average solution time for the uncorrelated instances was 516 seconds, and for the correlated ones, 380.

\subsection{Computational comparison among outsourcing policies}\label{subsect:cross-comparison}

The results discussed in the previous section do not allow a deep comparison between the different outsourcing policies being studied.
In this section we focus on such a comparison. We start by analyzing the termination status when using the different formulations.
Afterwards, we report percentage gaps at termination for instances that could not be solved to proven optimality and computing times for those that could.
Finally, we analyze the capability of each outsourcing policy to provide a good approximation for other policy(ies).
This is accomplished by considering the optimal location-allocation solution obtained using one policy and checking how good it is if some other policy is considered instead.

\begin{figure}[h!]
	\centering
	\subfigure[Facility outsourcing (FO).]{
		\includegraphics[clip,width=0.30\textwidth]{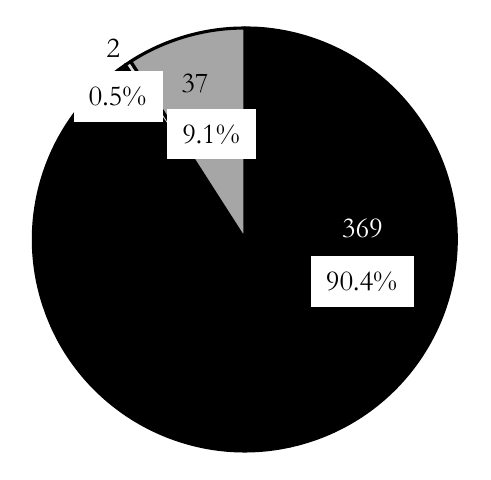}
		\label{terminationStatus-FO}
	}\quad
	\subfigure[Reassignment outsourcing (RO).]{
		\includegraphics[clip,width=0.30\textwidth]{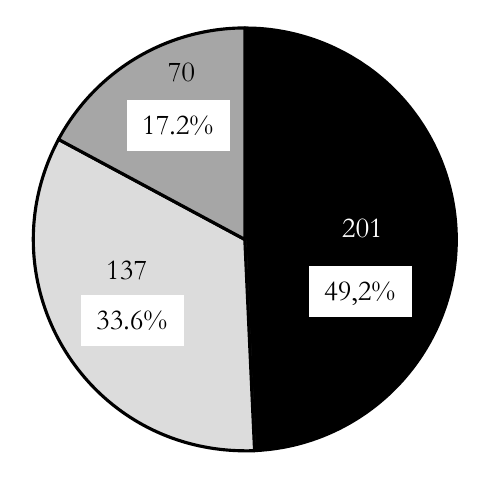}
		\label{terminationStatus-RO}
	} \\
	\subfigure[Order-driven customer outsourcing (OD-CO).]{
		\includegraphics[clip,width=0.30\textwidth]{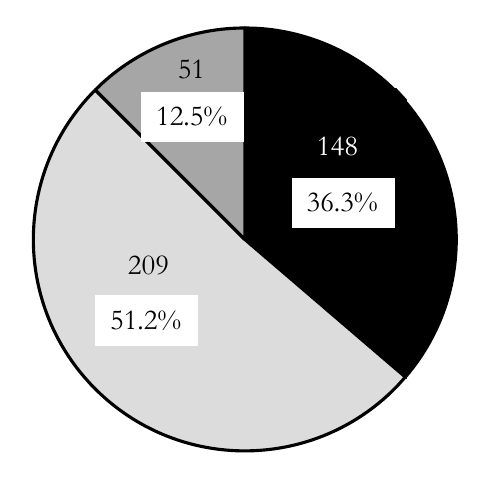}
		\label{terminationStatus-FIFO-CO}
	}\quad
	\subfigure[Cost-driven customer outsourcing (CD-CO).]{
		\includegraphics[clip,width=0.30\textwidth]{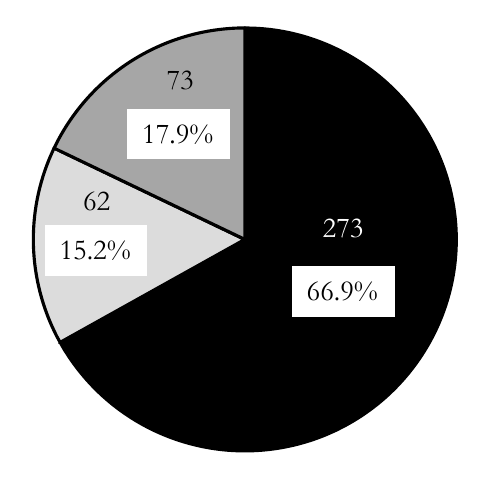}
		\label{terminationStatus-CD-CO}
	} \\
	\subfigure{
		\includegraphics[clip,width=0.50\textwidth]{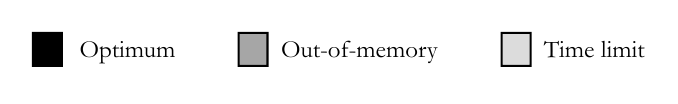}
	}
	\caption{Termination status for the different outsourcing strategies.}
	\label{figure:terminationStatus}
\end{figure}

Figure~\ref{figure:terminationStatus} depicts the termination status for the different outsourcing policies.
We can observe that the policy for which we were able to obtain more optimal solutions (around 90\%) is FO (Figure~\ref{terminationStatus-FO}).
In the other extreme, we observe OD-CO where optimality was proven for nearly 37\% of the instances (Figure~\ref{terminationStatus-FIFO-CO}).
When we compare this policy with the other CO policy, we conclude that the inclusion of the ordering constraints  \eqref{P2:order} in formulation \eqref{P2:of}--\eqref{P2:sbin} has a dramatic consequence in terms of model solvability.
In the case of RO and OD-CO, the major cause for not solving an instance to proven optimality is the time limit, whereas for the other two policies available memory was the reason for a premature stop. Of course, these numbers apply to the computer we have used in our experiments. Still, while memory limitation could be alleviated by using a computer with a higher RAM, the computing time limitation could also be alleviated by increasing the time limit. Thus, in our opinion, the above analysis reflects the main difficulties for optimally solve the instances for each of the considered outsourcing policies, and does not depend on the computer we have used.

\begin{figure}[h!]
	\centering
	\subfigure[Facility outsourcing (FO).]{
		\includegraphics[clip,width=0.45\textwidth]{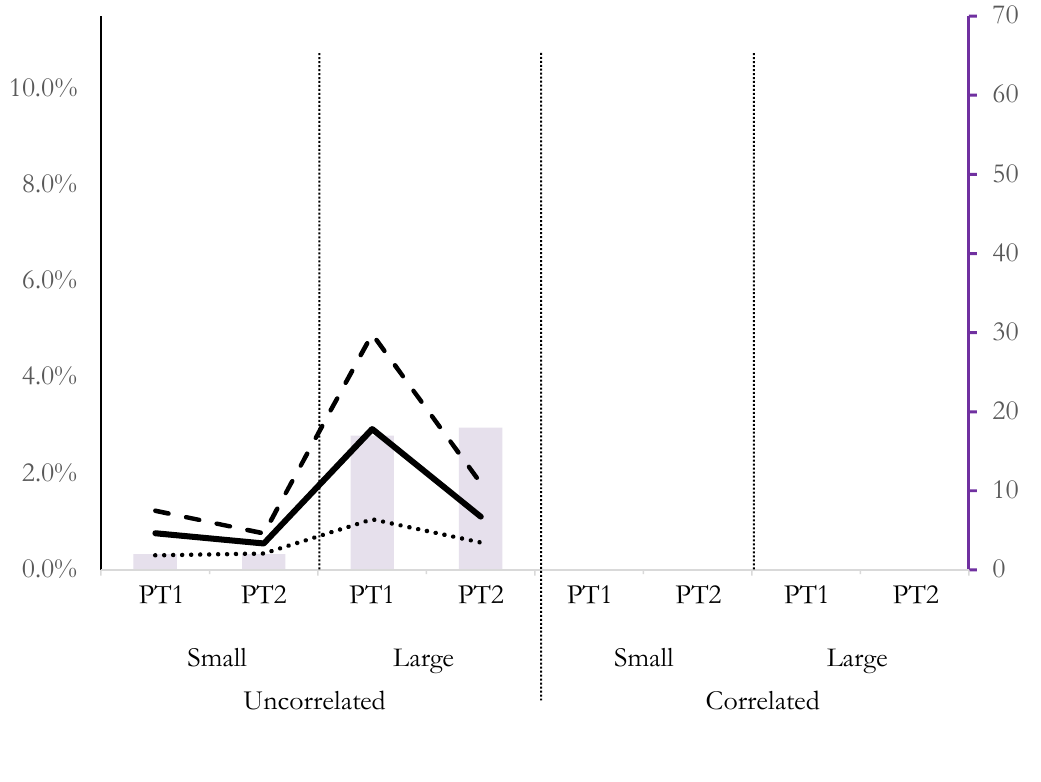}
		\label{terminationGap-FO}
	}\quad
	\subfigure[Reassignment outsourcing (RO).]{
		\includegraphics[clip,width=0.45\textwidth]{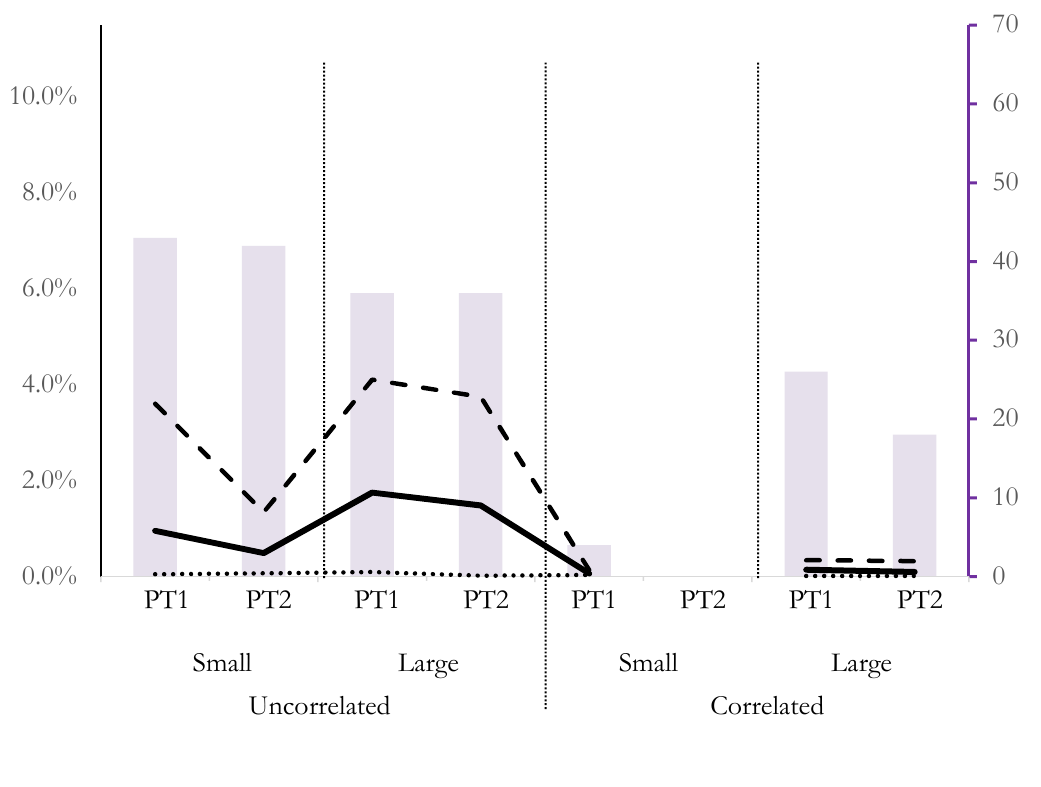}
		\label{terminationGap-RO}
	} \\
	\subfigure[Order-driven customer outsourcing (OD-CO).]{
		\includegraphics[clip,width=0.45\textwidth]{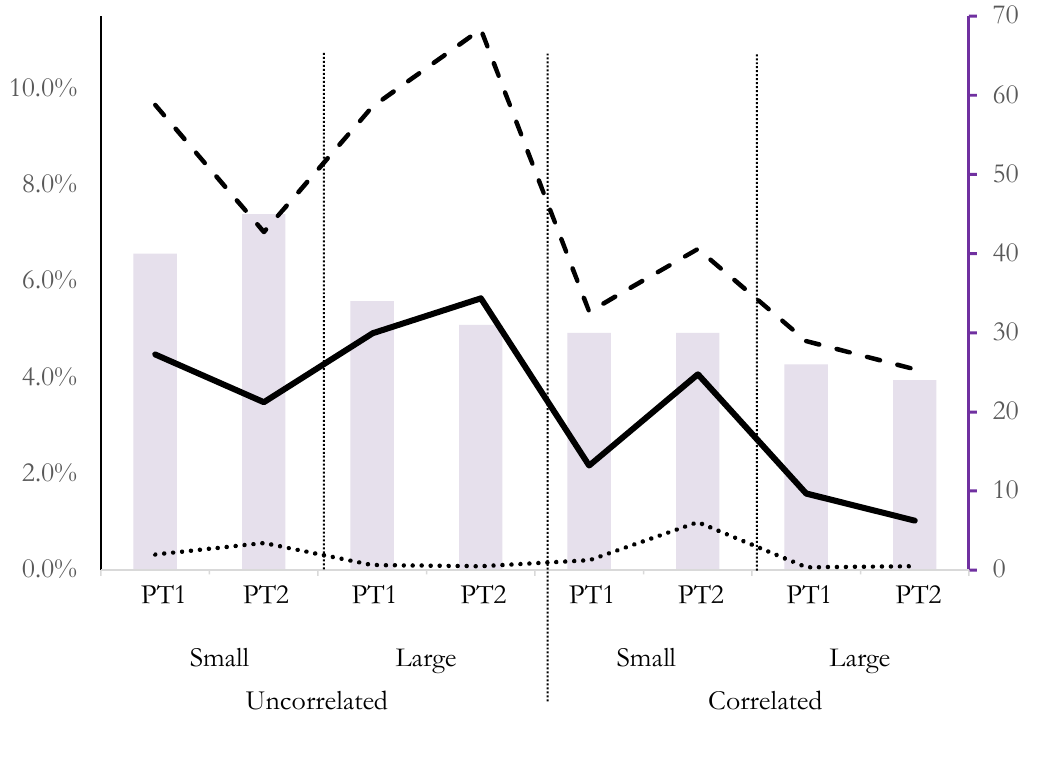}
		\label{terminationGap-OD-CO}
	}\quad
	\subfigure[Cost-driven-customer outsourcing (CD-CO).]{
		\includegraphics[clip,width=0.45\textwidth]{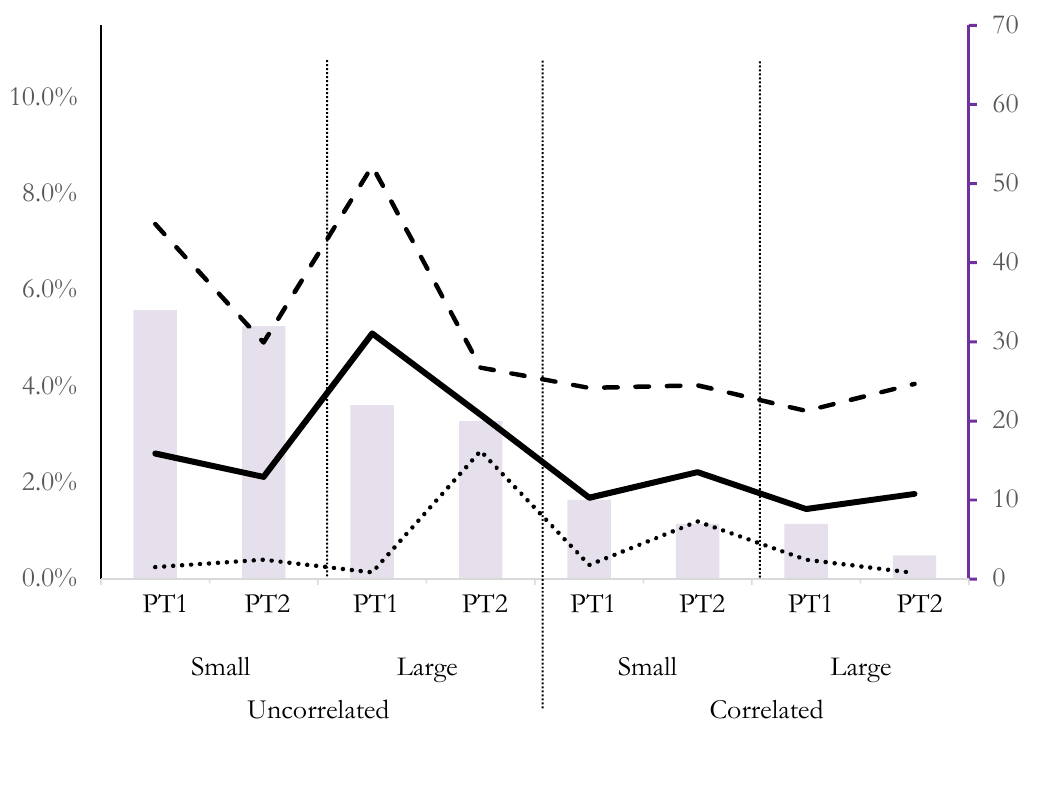}
		\label{terminationGap-CD-CO}
	}
	\caption{Percentage gaps at termination for the different outsourcing strategies.}
	\label{figure:terminationGap}
\end{figure}

In Figure~\ref{figure:terminationGap} we can observe the termination percentage gap (\%) for the instances that could not be solved to optimality either because the time limit was reached or because the computer memory was exhausted.
Each sub-figure refers to one outsourcing policy as labeled.
The bars indicate the number of instances involved (to be read in the right-hand side axis); The three lines (whose values can be seen in the left-hand side axis) indicate the minimum, average and maximum values observed.
In these figures, we are disaggretating the results according to three factors: type of demand (correlated or uncorrelated) dimension of the instances (small or large), and demand pattern (pattern 1 or pattern 2).
In these figures we do not distinguish between the two capacity types, low and high ($\gamma=1$ and $\gamma=2$) because for the purpose of the analysis of the percentage gaps, that distinction was not as relevant as the others when we look into the results.

Observing Figure~\ref{figure:terminationGap} we realize that the instances with correlated demand seem easier to handle when compared to those with uncorrelated demand. In fact, for the former, not only fewer instances are involved in these figures (i.e., more instances were solved to optimality) but also for those considered in the figure, the final gap at termination is globally smaller. As before, this is an indication that by working with correlated data one is considering an easier setting in terms of possibilities for the future observations since the customers behave in a similar fashion.

When we focus on the demand pattern or on the size of the instances, we do not observe a clear trend in Figure~\ref{figure:terminationGap}. This is an indication that the structure of the problem is not significantly influenced by these aspects. Overall, the final gap at termination for the instances not solved to optimality does not seem to be too sensitive neither to the dimension of the instances nor to the demand pattern.

When comparing the different outsourcing policies, we see a clear trend in favor of FO (Figure~\ref{terminationGap-FO}) and RO (Figure~\ref{terminationGap-RO}).
On the other hand, OD-CO seems again to provide the model harder to solve to proven optimality using a general purpose solver.
As mentioned before, the need to include ordering constraints in the model seems to make a huge difference in its structure with a dramatic consequence in terms of its efficient solvability.

\begin{figure}[h!]
	\centering
	\subfigure[Facility outsourcing (FO).]{
		\includegraphics[clip,width=0.45\textwidth]{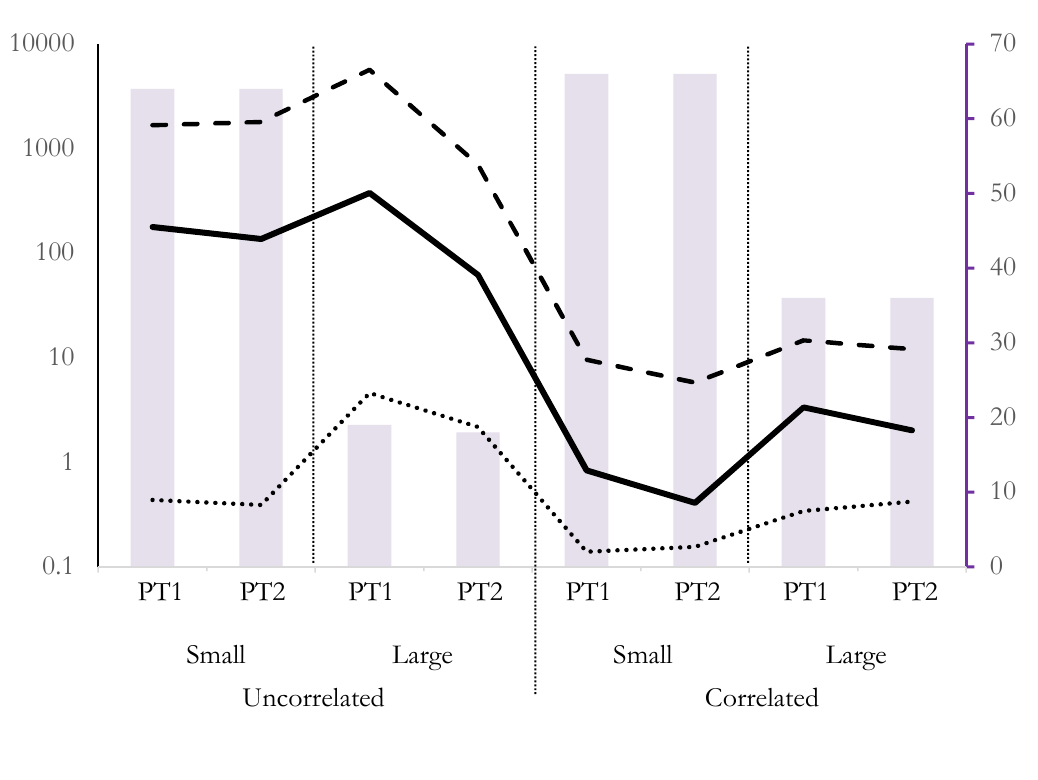}
		\label{terminationCPU-FO}
	}\quad
	\subfigure[Reassignment outsourcing (RO).]{
		\includegraphics[clip,width=0.45\textwidth]{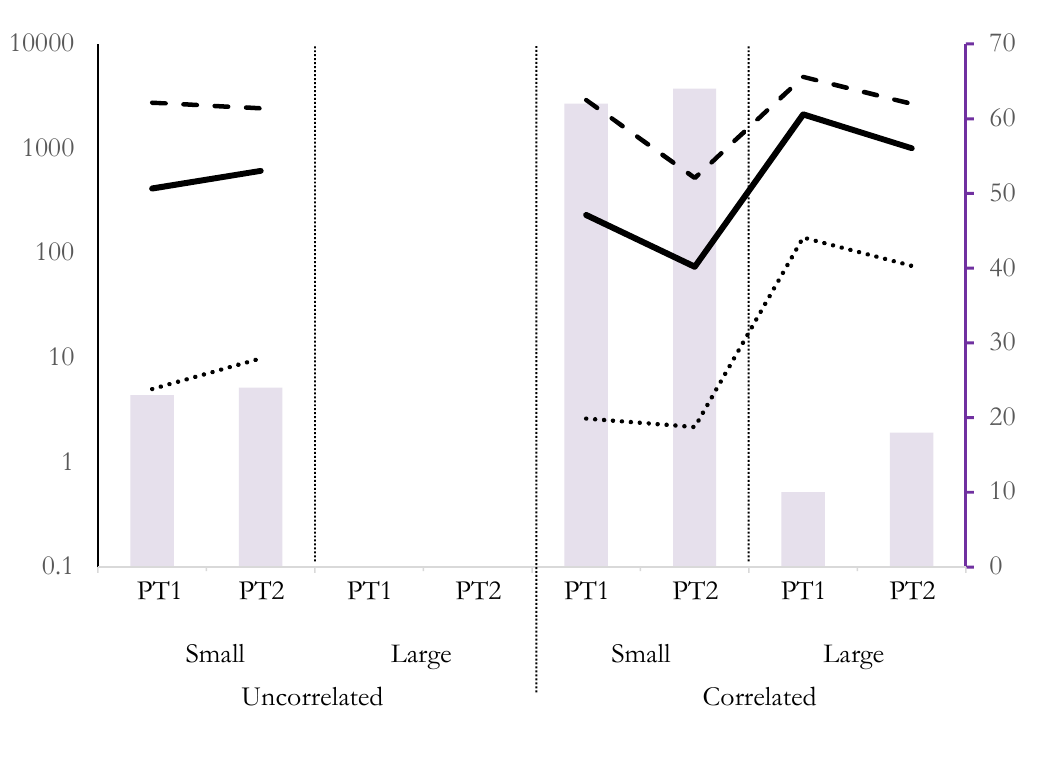}
		\label{terminationCPU-R}
	} \\
	\subfigure[Order-driven customer outsourcing (OD-CO).]{
		\includegraphics[clip,width=0.45\textwidth]{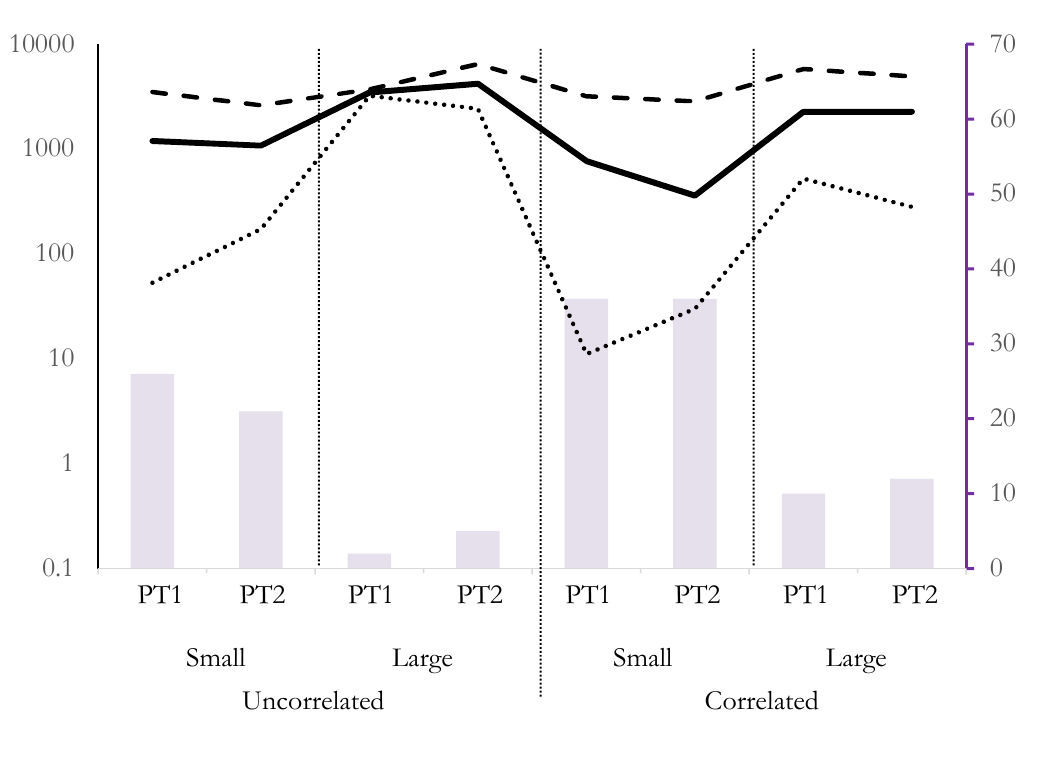}
		\label{terminationCPU-OD-CO}
	}\quad
	\subfigure[cost-driven customer outsourcing (CD-CO).]{
		\includegraphics[clip,width=0.45\textwidth]{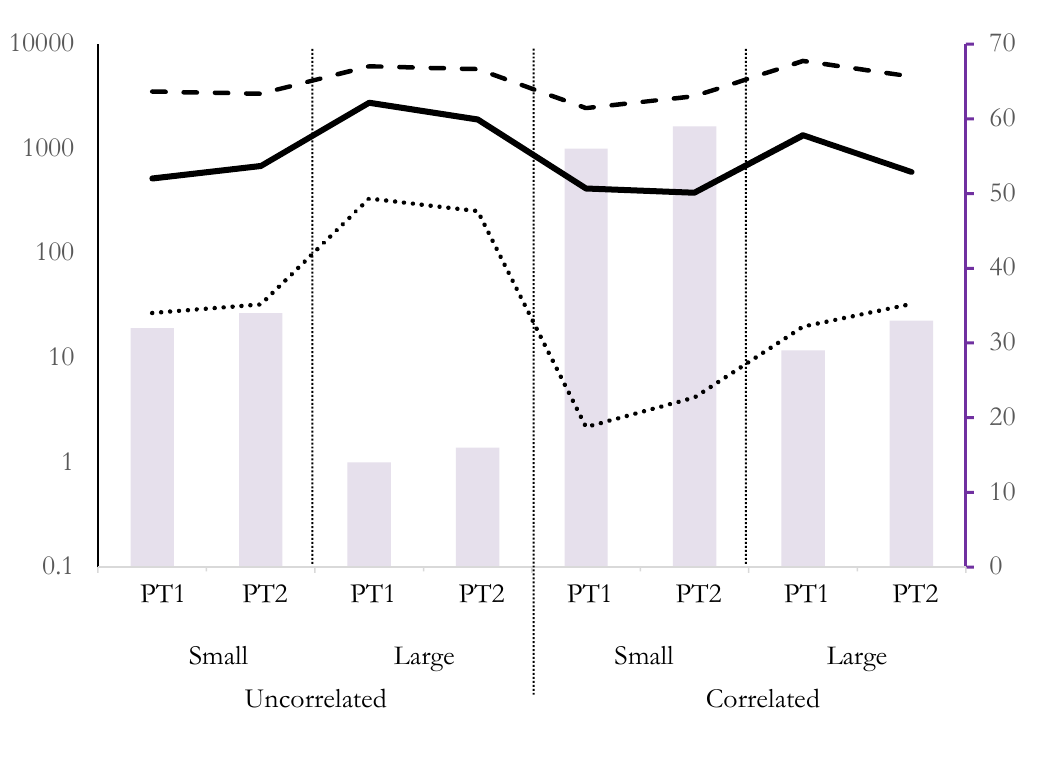}
		\label{terminationCPU-CD-CO}
	}
	\caption{Computing times (in seconds) for the instances that were successfully solved to proven optimality.}
	\label{figure:terminationCPU}
\end{figure}

In Figure~\ref{figure:terminationCPU} we can observe the computing times (in seconds) for the instances that were successfully solved up to proven optimality  within the time limit. Each sub-figure refers to one outsourcing policy as labeled.
The vertical bars count the number of instances involved in the analysis whose corresponding values are marked on the right-hand side axis.
The three lines depicted (whose values are associated with the left-hand side axis) regard the minimum, average and maximum values observed.
Similarly to Figure~\ref{figure:terminationGap} we are disaggregating the results according to type of demand (uncorrelated or correlated), size of the instances and demand pattern.
As before, FO is the policy whose optimization model is easier to tackle.
This is even more evident when we focus on the instances with correlated data.
Apart from FO, the type of demand does not seem to influence the ``friendliness'' of the optimization models.
Analyzing the graphics we also see a tendency for fewer instances involved in the analysis for large instances with uncorrelated demand, which is in line to previous conclusions already drawn.

\begin{figure}[p]
	\centering
	\subfigure[FO, uncorrelated data.]{
		\includegraphics[clip,width=0.60\textwidth]{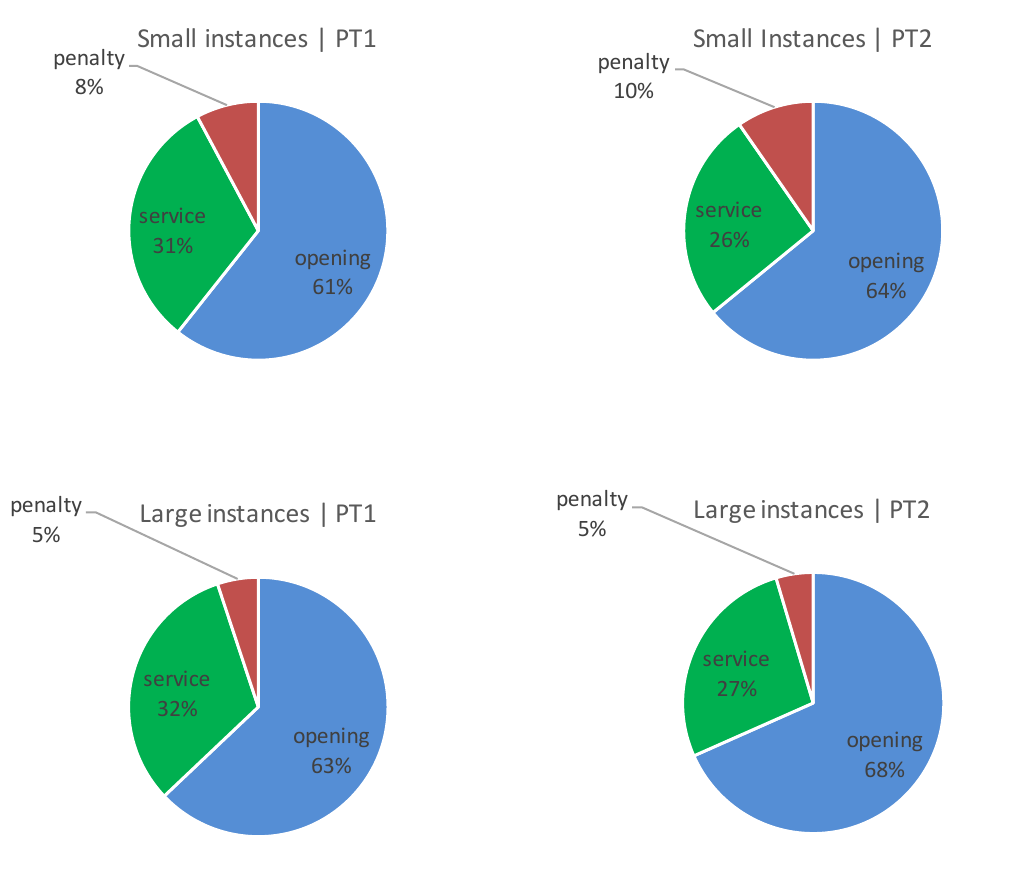}
		\label{CostStructure-FO-U}
	} \\
	\subfigure[FO, correlated data..]{
		\includegraphics[clip,width=0.60\textwidth]{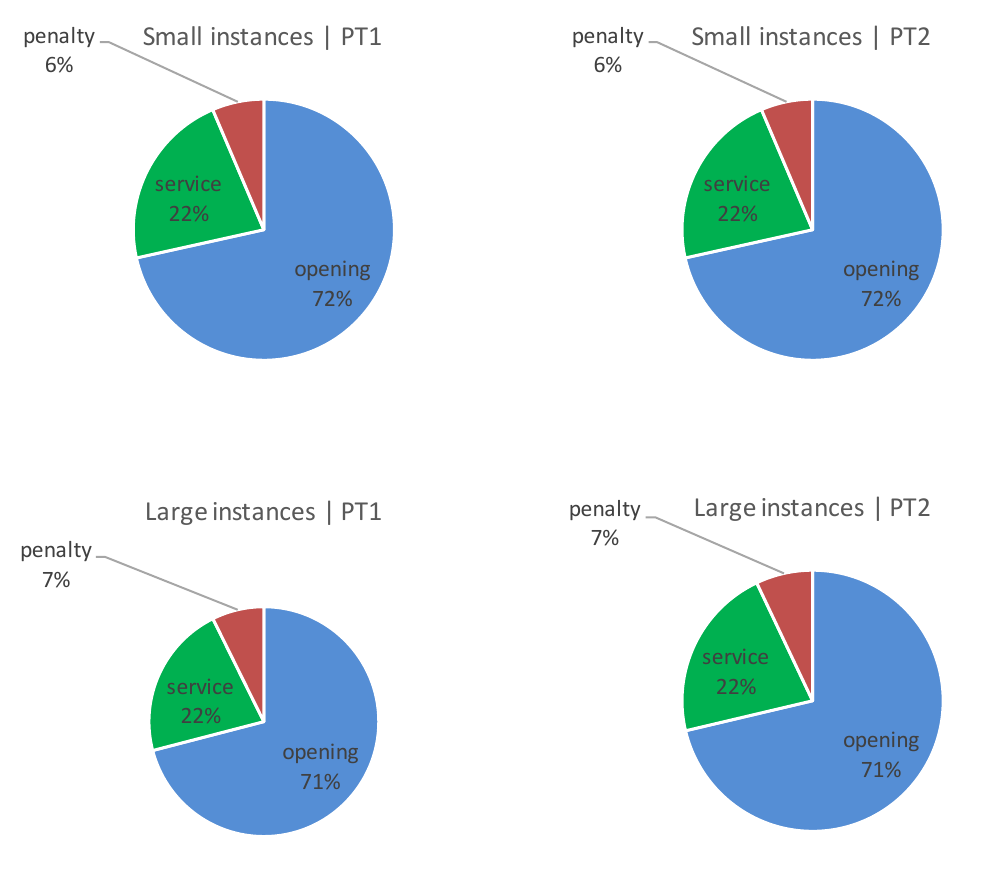}
		\label{CostStructure-FO-C}
	}
	\caption{Cost structure---FO.}
	\label{CostStructure-FO}
\end{figure}

\begin{figure}[p]
	\centering
	\subfigure[OD-CO, uncorrelated data.]{
		\includegraphics[clip,width=0.60\textwidth]{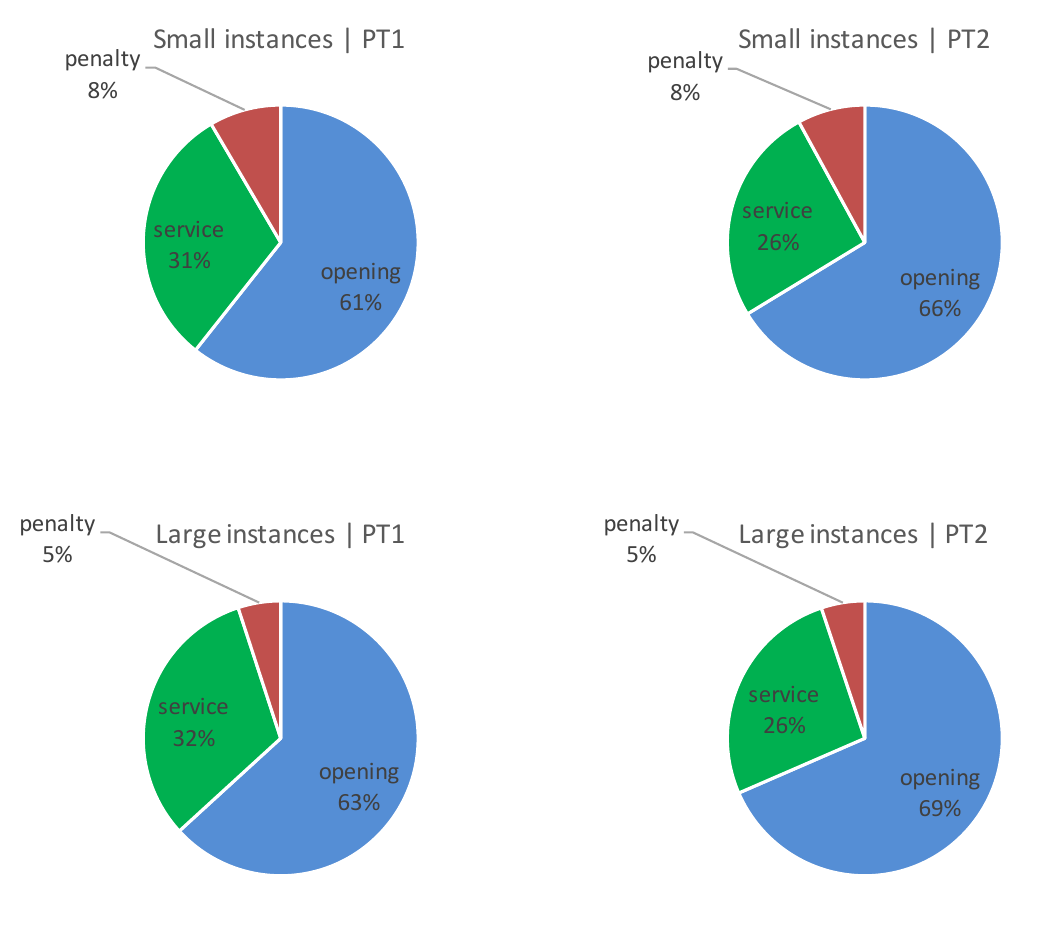}
		\label{CostStructure-OD-CO-U}
	} \\
	\subfigure[OD-CO, correlated data.]{
		\includegraphics[clip,width=0.60\textwidth]{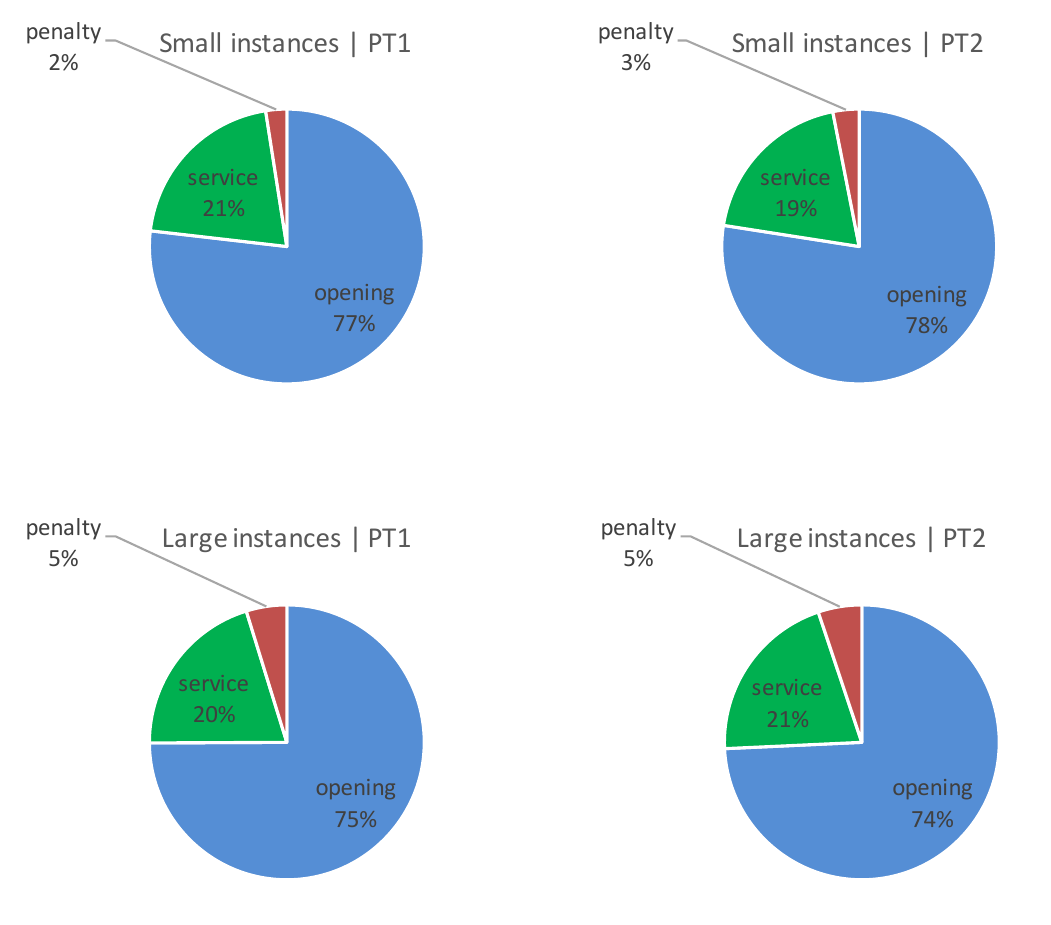}
		\label{CostStructure-OD-CO-C}
	}
	\caption{Cost structure---OD-CO.}
	\label{CostStructure-OD-CO}
\end{figure}

\begin{figure}[p]
	\centering
	\subfigure[CD-CO, uncorrelated data.]{
		\includegraphics[clip,width=0.60\textwidth]{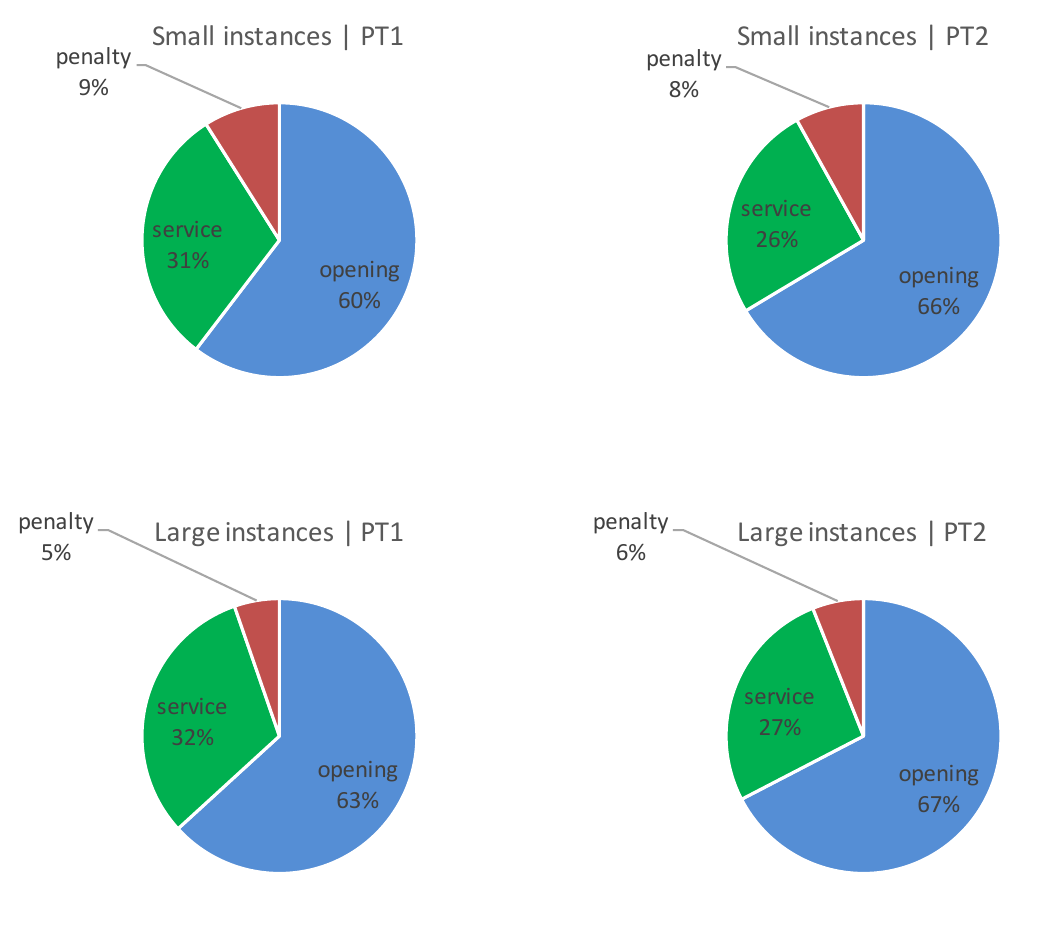}
		\label{CostStructure-CD-CO-U}
	} \\
	\subfigure[CD-CO, correlated data.]{
		\includegraphics[clip,width=0.60\textwidth]{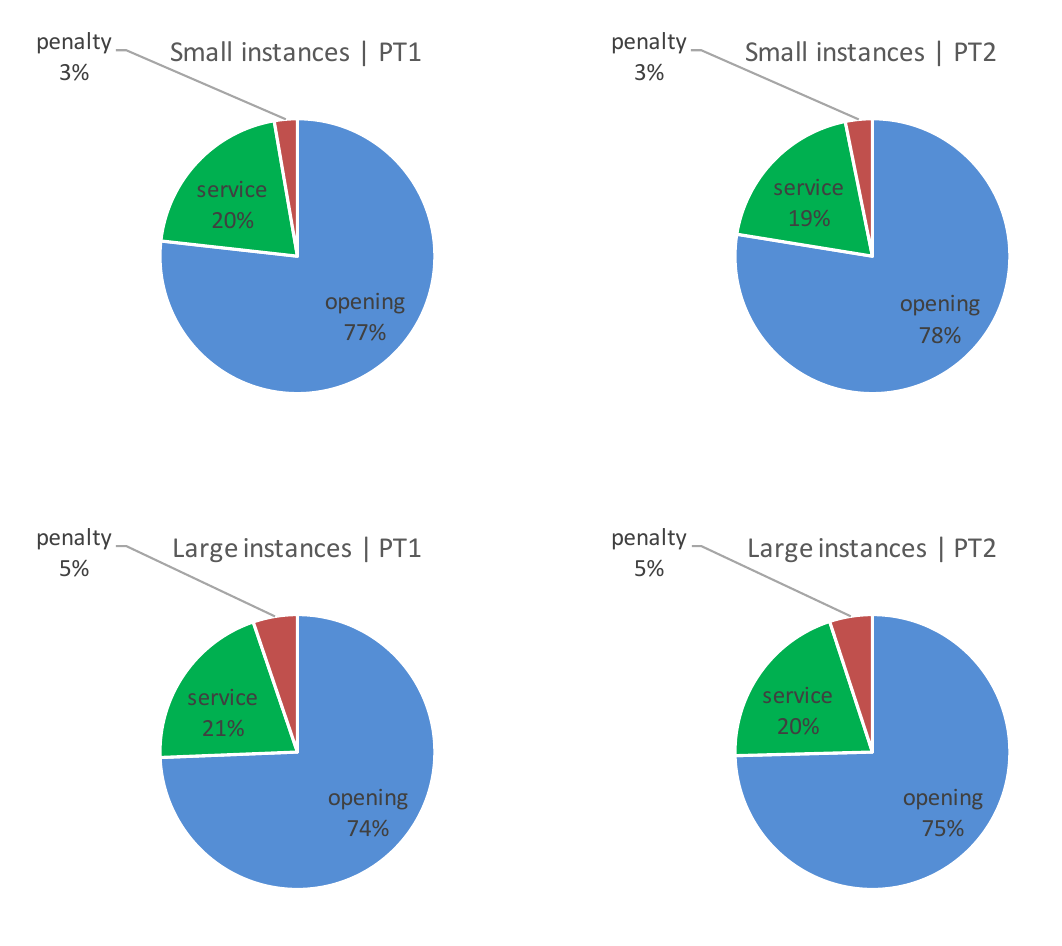}
		\label{CostStructure-CD-CO-C}
	}
	\caption{Cost structure---CD-CO.}
	\label{CostStructure-CD-CO}
\end{figure}

\begin{figure}[p]
	\centering
	\subfigure[RO, uncorrelated data.]{
		\includegraphics[clip,width=0.65\textwidth]{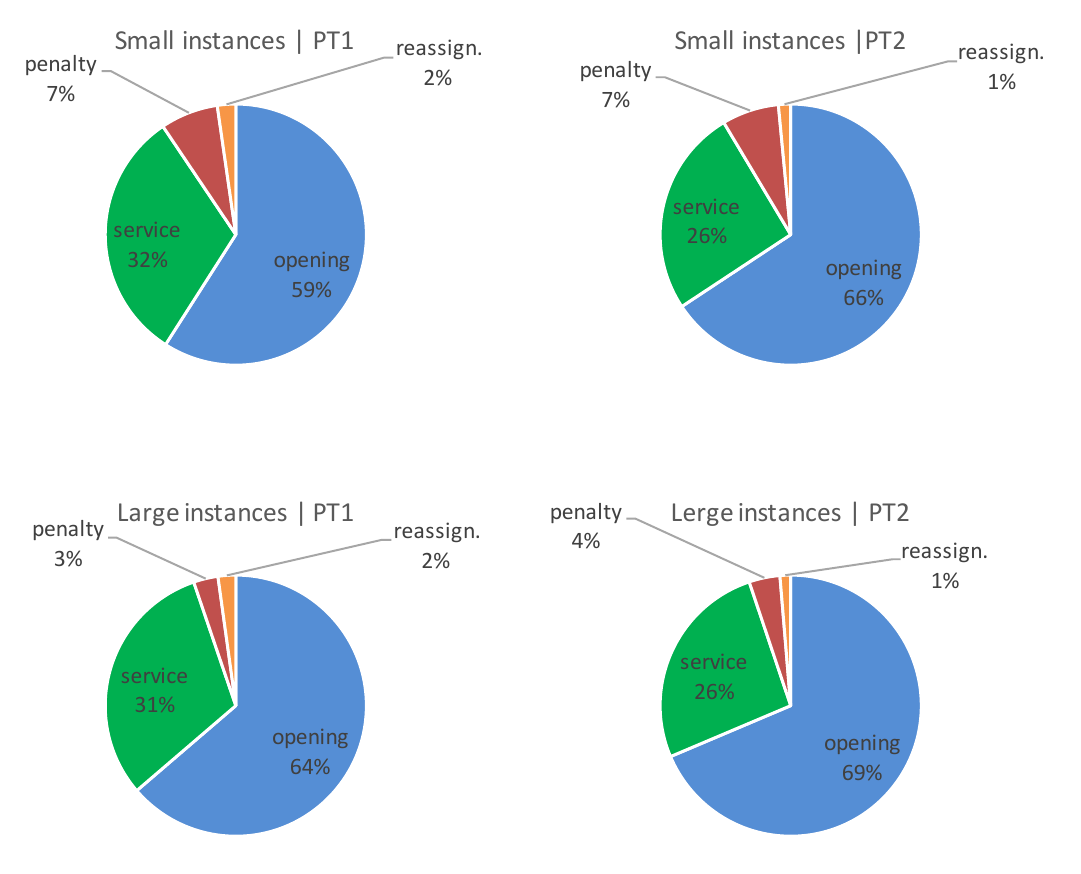}
		\label{CostStructure-RO-U}
	} \\
	\subfigure[RO, correlated data.]{
		\includegraphics[clip,width=0.65\textwidth]{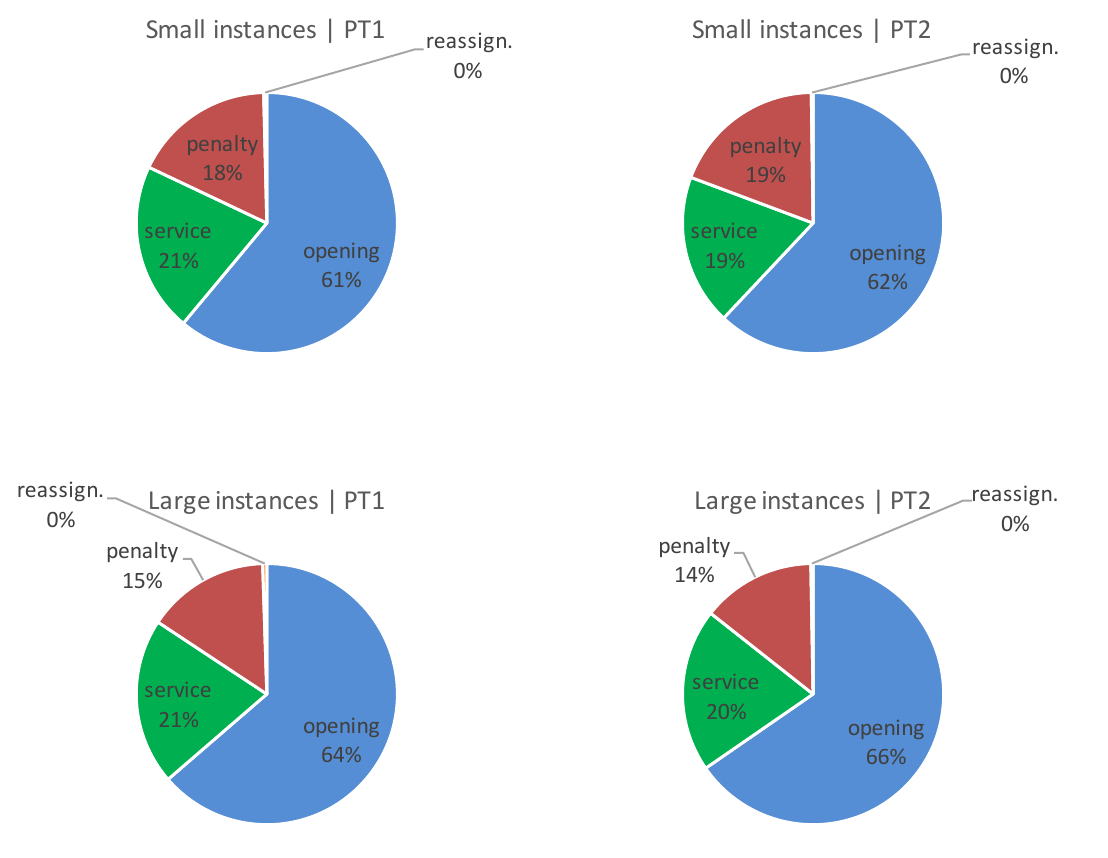}
		\label{CostStructure-RO-C}
	}
	\caption{Cost structure---RO.}
	\label{CostStructure-RO}
\end{figure}

Finally, in this section we look into the cost structure of the optimal solutions found using each outsourcing strategy.
In particular, for each outsourcing policy we present the contribution to the objective function value of the strategic decisions (facility setup costs), transportation costs (customers assignment costs) and outsourcing costs.
This information is depicted in Figures~\ref{CostStructure-FO}--\ref{CostStructure-RO}.
In these figures ``opening'' refers to the strategic decisions, ``service'' concerns customer satisfaction cost, and ``penalty'' stands for the outsourcing.
In the case of facility reassignment, we also present the corresponding costs---``reassign.''.

Observing Figures~\ref{CostStructure-FO}--\ref{CostStructure-RO} we immediately conclude for a tendency in terms of higher opening costs when data is correlated.
This is an indication that under correlated data, we need to pay more in terms of first-stage decisions (facility opening costs) to better hedge against the future uncertainty.
On the other hand, apart from RO, we see that the penalty costs are usually small and even more with correlated data.
Regarding RO we do not see much different in terms of the opening costs when moving from uncorrelated to correlated data.
However, in the latter the reassignment costs are negligible which is an indication that under correlated data reassigning customers is less frequent, which is justified by having most of the facilities close to their service capacity.

\subsection{Managerial insight}
The results reported so far in this section show that the mathematical programming formulations adopted for the different outsourcing policies have a different behavior when it comes to tackling them using a general purpose solver.
Overall, there seems to be a hierarchy between the four models tested.
FO seems to provide the model easier to solve optimally by the solver; then we observe RO; finally, we see the customer outsourcing policies with OD-CO associated with the model harder to solve optimally.\\

Given the differences observed between the optimization models, a decision maker may think of using the ones easier to solve to proven optimally to provide approximate solutions to the harder ones. To investigate this possibility, for every instance we looked into how good the optimal (or best) solution found for one outsourcing policy is when looked at as an approximate solution for another policy.
In Table~\ref{table:crossed_results-global} we can observe the values obtained.

\begin{table}[!h]
	\begin{center}
		{\footnotesize
			\begin{tabular}{l c r c r c r c r }
				\toprule
				& &     FO & & CD-CO & &  OD-CO & &     RO \\
				\cmidrule{3-9}
				FO & &  \gray & & 0.463 & & 20.378 & & -0.014 \\
				CD-CO & & 17.826 & & \gray & & 48.583 & & 13.169 \\
				OD-CO & &  2.896 & & 1.646 & &  \gray & &  5.318 \\
				RO & &  1.139 & & 3.044 & & 23.039 & &  \gray \\
				\bottomrule
			\end{tabular}
		}
	\end{center}
	\caption{Average percentage gaps (\%) of the optimal/best-found solution using one policy (rows) with respect to the others (columns).}\label{table:crossed_results-global}
\end{table}

Looking at the values per rows in this table, we conclude that the policy that provides more robust solutions for all the other policies is OD-CO.
In fact, apart from this one, all the other policies provide a quite bad approximation for at least one other outsourcing policy.
If we focus on the values per columns, we conclude that CD-CO is the policy with the optimal solution ``easier'' to approximate using the solution provided by other policies. An interesting aspect to notice is that the ``matrix'' presented in Table~\ref{table:crossed_results-global} is not symmetric.
The differences are often quite significant.
This indicates that on average, the fact that one policy provides a good approximation to another one does not imply the reverse.
The most extreme case involves the two customer outsourcing policies: we see that OD-CO solutions are very good approximations for CD-CO but the reverse approximation is quite bad.
This indicates that the optimal (or best feasible solutions) are clearly more sensitive to order-driven outsourcing than to cost-driven outsourcing.
We can draw similar conclusions for other pairs of policies.

In order to deepen the analysis, we disaggregated the results.
In Table~\ref{table:crossed_results-K} we consider explicitly the type of demand (uncorrelated or correlated) and the type of capacity adopted for the facilities---low or high ($\gamma=1$, $\gamma=4$).
In Table~\ref{table:crossed_results-p} we consider again the type of demand but now jointly with the demand pattern as well.
We note that in many cases, there were instances that could not be solved to proven optimality (either because the time limit was reached or because the memory of the machine was exhausted).
In such situations we work with best-found solutions.
This explains some negative averages that we observe in Tables~\ref{table:crossed_results-K} and \ref{table:crossed_results-p}.

Observing Tables~\ref{table:crossed_results-K} and \ref{table:crossed_results-p} we conclude that in general, uncorrelated demand leads to worse results than correlated demand.
In other words, the different outsourcing policies provide poorer approximations to each other in the former case.
In these tables we also observe that changing the capacity of the facilities seems to have a greater impact than changing the demand pattern.
Additionally, we realize that the results for $\gamma=4$ (high capacities) seem systematically worse than for $\gamma=1$ (low capacities).
Moreover, when we focus on OD-CO we observe that the instances with uncorrelated demand seem invariably more difficult to approximate using the other outsourcing policies than those with correlated demand.
The worst approximation occurs (on average) when we approximate OD-CO by CD-CO for uncorrelated demand and the larger capacity of the facilities.
The best approximations (on average) occur when we approximate RO by FO.

\begin{table}[h!]
	\begin{center}
		{\footnotesize
			\begin{tabular}{ll c rr c rr c rr c rr }
				\toprule
				&              & & \multicolumn{2}{c}{FO} & & \multicolumn{2}{c}{CD-CO} & & \multicolumn{2}{c}{OD-CO} & & \multicolumn{2}{c}{RO} \\
				\cmidrule{4-5} \cmidrule{7-8} \cmidrule{10-11} \cmidrule{13-14}
				&              & &           $\gamma$=1 &    $\gamma$=4 & &               $\gamma$=1 &   $\gamma$=4 & &             $\gamma$=1 &     $\gamma$=4 & &         $\gamma$=1 &     $\gamma$=4 \\
				\cmidrule{4-5} \cmidrule{7-8} \cmidrule{10-11} \cmidrule{13-14}
				\multirow{2}{*}{FO}
				& U & &         \gray &  \gray & &             0.335 & 0.012 & &          17.526 &  28.233 & &      -0.058 &   0.000 \\
				&  C & &         \gray &  \gray & &             1.506 & 0.000 & &          10.672 &  25.081 & &       0.000 &   0.000 \\
				&              & &               &        & &                   &       & &                 &         & &             &         \\
				\multirow{2}{*}{CD-CO}
				& U & &         4.708 & 52.415 & &             \gray & \gray & &          28.325 & 110.513 & &       2.962 &  40.207 \\
				& C & &         2.854 & 11.327 & &             \gray & \gray & &          12.080 &  43.415 & &       3.581 &   5.924 \\
				&              & &               &        & &                   &       & &                 &         & &             &         \\
				\multirow{2}{*}{OD-CO}
				& U & &         0.893 &  0.230 & &             0.005 & 0.121 & &           \gray &   \gray & &       3.786 &   1.089 \\
				& C & &         9.854 &  0.608 & &             6.045 & 0.411 & &           \gray &   \gray & &      14.267 &   2.131 \\
				&              & &               &        & &                   &       & &                 &         & &             &         \\
				\multirow{2}{*}{RO}
				& U & &         1.316 &  0.160 & &             3.600 & 0.160 & &          21.210 &  28.366 & &       \gray &   \gray \\
				& C & &         1.242 &  1.837 & &             4.720 & 3.700 & &          13.527 &  29.051 & &       \gray &   \gray \\
				\bottomrule
			\end{tabular}
		}
	\end{center}
	\caption{Data of Table~\ref{table:crossed_results-global} disaggregated according to the type of data and the type of capacity for the facilities.}\label{table:crossed_results-K}
\end{table}

\begin{table}[h!]
	\begin{center}
		{\footnotesize
			\begin{tabular}{ll c rr c rr c rr c rr }
				\toprule
				&              & & \multicolumn{2}{c}{FO} & & \multicolumn{2}{c}{CD-CO} & & \multicolumn{2}{c}{OD-CO} & & \multicolumn{2}{c}{RO} \\
				\cmidrule{4-5} \cmidrule{7-8} \cmidrule{10-11} \cmidrule{13-14}
				&              & &           PT1 &    PT2 & &              PT1 &   PT2 & &             PT1 &     PT2 & &         PT1 &    PT2  \\
				\cmidrule{4-5} \cmidrule{7-8} \cmidrule{10-11} \cmidrule{13-14}
				\multirow{2}{*}{FO}
				& U & &         \gray &  \gray & &            -0.145 & 0.491 & &          22.291 &  23.468 & &      -0.027 &  -0.031 \\
				& C & &         \gray &  \gray & &             0.740 & 0.766 & &          16.714 &  19.039 & &       0.000 &   0.000 \\
				&              & &               &        & &                   &       & &                 &         & &             &         \\
				\multirow{2}{*}{CD-CO}
				& U & &        28.150 & 28.973 & &             \gray & \gray & &          73.127 &  65.711 & &      21.653 &  21.516 \\
				& C & &         7.336 &  6.845 & &             \gray & \gray & &          29.862 &  25.633 & &       5.275 &   4.231 \\
				&              & &               &        & &                   &       & &                 &         & &             &         \\
				\multirow{2}{*}{OD-CO}
				& U & &         0.546 &  0.576 & &             0.196 &-0.070 & &           \gray &   \gray & &       3.786 &   2.093 \\
				& C & &         3.797 &  6.664 & &             2.310 & 4.147 & &           \gray &   \gray & &       7.565 &   8.833 \\
				&              & &               &        & &                   &       & &                 &         & &             &         \\
				\multirow{2}{*}{RO}
				& U & &         1.316 &  0.363 & &             3.601 & 1.166 & &          21.210 &  24.262 & &       \gray &   \gray \\
				& C & &         1.137 &  1.942 & &             3.887 & 4.525 & &          20.372 &  22.205 & &       \gray &   \gray \\
				\bottomrule
			\end{tabular}
		}
	\end{center}
	\caption{Data of Table~\ref{table:crossed_results-global} disaggregated according to the type of data and the demand pattern.}\label{table:crossed_results-p}
\end{table}

\begin{figure}[!h]
	\centering
	\subfigure[Disaggregation according to type of data, size of instance, and capacity of the facilities.]{
		\includegraphics[clip,width=0.65\textwidth]{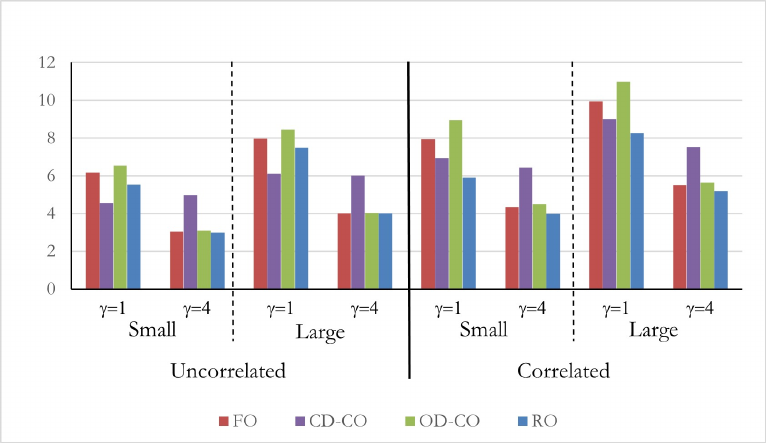}
		\label{FacilitiesOpen_K}
	}\\
	\subfigure[Disaggregation according to type of data, size of instance, and demand pattern.]{
		\includegraphics[clip,width=0.65\textwidth]{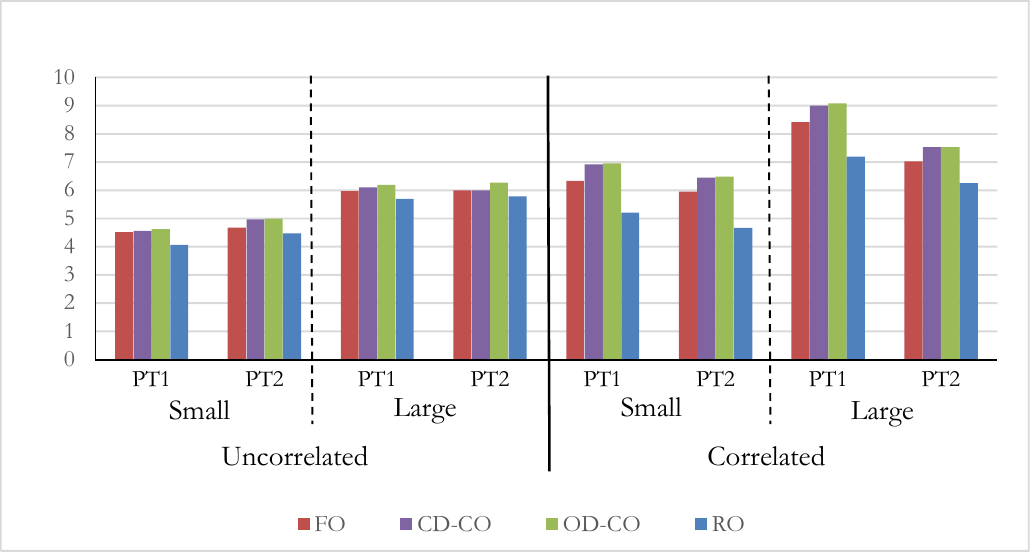}
		\label{FacilitiesOpen_PT}
	}
	\caption{Average number of open facilities in the optimal or best feasible solution found.}
	\label{figure:AvgNumberFacilities}
\end{figure}

In search for an explanation in terms of the robustness of model FO (for providing good feasible solutions for CD-CO and RO) and OD-CO for all policies, we analyzed the number of facilities opened in each solution.
Figure~\ref{figure:AvgNumberFacilities} depicts such results already disaggregated according to the type of demand (uncorrelated or correlated) and the size of the instances.
Additionally, in Figure~\ref{FacilitiesOpen_K} we can observe the results for the two types of capacity whereas in  Figure~\ref{FacilitiesOpen_PT} we have the information according to the demand pattern.

Observing this figure we conclude that OD-CO policy calls (on average) for the largest number of open facilities.
This may turn the solutions for this policy more robust when considered as feasible solutions to the other policies since more facilities means, for instance, more possibilities for reassigning customers and lower customer outsourcing costs.

\section{Conclusions} \label{sec:conclu}

In this work we studied different outsourcing policies in the context of the facility location problem with Bernoulli demands.
We extended the previous work on this problem by considering both uncorrelated and correlated service demands.
Furthermore, in addition to the two policies that had already been proposed in the literature we introduced two other possibilities.
Each of the four investigated outsourcing policies leads to a variant of the problem for which a mathematical programming formulation was derived assuming a finite set of scenarios for the demand.
An extensive computational study was performed to evaluate the extent to which the mixed-integer linear programming models proposed can be handled by a general-purpose solver as well as to evaluate the capability of each model to produce high quality feasible solutions for the others.

The results show that two outsourcing policies lead to models easier to solve to proven optimality using a solver, namely: facility outsourcing and reassignment outsourcing. Moreover, the fact that the two customer outsourcing policies induce mathematical models more difficult to solve to proven optimality, motivates the use of approximate solutions in those cases. The results show that on average, the optimal solution obtained for facility outsourcing turns out to be a good approximate solution for cost-driven customer outsourcing and for reassignment outsourcing.
Nevertheless, the most ``robust'' outsourcing policy is order-driven customer outsourcing in the sense that the corresponding optimal solutions are invariably good feasible solutions to the other three policies. This means that if a decision maker adopts the a priori solution provided by an order-driven customer outsourcing, it is likely to perform reasonably well, regardless the outsourcing policy that is finally adopted.

This work opens new research directions in terms of the explicit inclusion of outsourcing in stochastic capacitated facility location problems.
In particular, it would be interesting to investigate whether probability distributions other than Bernoulli can benefit from the insights provided by this paper.

\section*{Acknowledgments}
This work was supported by the Spanish Ministry of Economy and Competitiveness through MINECO/ FEDER grants MTM2015-63779-R and MTM2019-105824GB-I00, and by National Funding from FCT---Fundação para a Ciência e a Tecnologia, Portugal, under the project: UIDB/04561/2020. This support is gratefully acknowledged.

The authors thank the two anonymous reviewers whose comments and insights helped improving the article.

\bibliography{references}

\end{document}